\documentclass{amsart}
\usepackage{amsfonts,latexsym,amsmath,amscd,amssymb}
\usepackage{bbm}

\theoremstyle{plain}
\newtheorem{theorem}{Theorem}
\numberwithin{theorem}{section}

\newtheorem{corollary}{Corollary}
\numberwithin{corollary}{section}

\newtheorem{definition}{Definition}
\numberwithin{definition}{section}

\newtheorem{lemma}{Lemma}
\numberwithin{lemma}{section}

\newtheorem{proposition}{Proposition}
\numberwithin{proposition}{section}

\newtheorem{remark}{Remark}
\numberwithin{remark}{section}

 \numberwithin{equation}{section}

\newcommand {\be}{\begin{equation}}
\newcommand {\ee}{\end{equation}}

\newcommand{\h}{\begin{eqnarray*}}
 \newcommand{\e}{\end{eqnarray*}}

\newcommand{\CC}{\mathbf{C}}
\begin{document}
\title[Modular Invariance and Twisted Anomaly Cancellations]{Modular Invariance and Twisted Anomaly Cancellations For Characteristic Numbers}

\author{Qingtao Chen}
\address{Q. Chen, Department of Mathematics, University of California,
Berkeley, CA, 94720-3840} \email{chenqtao@math.berkeley.edu}
\date{May 10, 2006}

\author{Fei Han}
\address{F. Han, \ Department of Mathematics, University of California,
Berkeley, CA, 94720-3840} \email{feihan@math.berkeley.edu}


\maketitle

\begin{abstract} By studying modular invariance properties of some characteristic forms, we obtain
twisted anomaly cancellation formulas. We apply these twisted
cancellation formulas to study divisibilities on spin manifolds
and congruences on spin$^c$ manifolds. Especially, we get twisted
Rokhlin congruences for $8k+4$ dimensional spin$^c$ manifolds.
\end{abstract}

\section {Introduction}

Let $M$ be a 12 dimensional smooth Riemannian manifold. A
beautiful relation between the top degree components of the
Hirzebruch $\widehat{L}$-form and $\widehat{A}$-form of $M$ was
shown by Alvarez-Gaum\'e and Witten [1] as a gravitational anomaly
cancellation formula as follows,\newline \be
\left\{\widehat{L}(TM,\nabla^{TM})\right\}^{(12)}=\left\{8\widehat{A}
(TM,\nabla^{TM})\mathrm{ch}(T_\mathbf{C}M,\nabla^{T_\mathbf{C}M})
-32\widehat{A}(TM,\nabla^{TM})\right\}^{(12)},\ee where
$T_\mathbf{C}M$ denotes the complexification of $TM$ and
$\nabla^{T_\mathbf{C}M}$ is canonically induced from
$\nabla^{TM}$, the Levi-Civita connection associated to the
Riemannian structure of $M$; $\mathrm{ch}(T_\mathbf{C}M,
\nabla^{T_\mathbf{C}M})$ denotes the Chern character form
associated to $(T_\mathbf{C}M, \nabla^{T_\mathbf{C}M})$ (cf.
[22]). This gravitational anomaly cancellation formula, which they
called miraculous cancellation formula, was derived from very
non-trivial computations.

(1.1) is generalized by Kefeng Liu [16] to arbitrary $8k+4$
dimensional manifolds by developing modular invariance properties
of characteristic forms. In [16], he proved that for each
$(8k+4)$-dimensional smooth Riemannian manifold $M$ the following
identity holds, \be \left \{\widehat{L}(TM, \nabla^{TM})\right
\}^{(8k+4)}=8\sum_{r=0}^{k}2^{6k-6r}h_r(T_\CC M).\ee In (1.2),
each $h_r(T_\CC M)$ is a differential form
$$\left \{\widehat{A}(TM, \nabla^{TM})\mathrm{ch}\left(b_r(T_\CC M), \nabla^{b_r(T_\CC M)}\right)\right\}^{(8k+4)},$$ where $b_r(T_\CC M)\in KO(M)\otimes \mathbf{C}, 0\leq r\leq k$, can
be derived canonically from $TM$. When the manifold is closed and
spin, according to the Atiyah-Hirzebruch divisibility [3],
$\langle h_r(T_\CC M), [M]\rangle$ are all even numbers. Therefore
from (1.2) and the Hirzebruch signature theorem [11],  we easily
get the Ochanine divisibility [19], which says that the signature
of an $8k+4$ dimensional smooth closed spin manifold is divisible
by 16. This shows us how miraculous cancellation formulas imply
divisibility of characteristic numbers. The author also provides a
similar miraculous cancellation formula for $8k$ dimensional
manifolds. Unfortunately this $8k$ dimensional cancellation
formula does not imply divisibility results. Note that in some
sense, Liu's formula refines the arguments of Hirzebruch [10] and
Landweber [14], who deduce the Ochanine divisibility by using the
ideas of elliptic genus, to the level of differential forms. See
[16] for details.

Liu's method is taken over in [8, 9] to study the Ochanine
congruence and the Finashin congruence. The authors show that
there actually exist more general miraculous cancellation formulas
with an extra complex line bundle involved, which turn out to be
efficient to study congruence phenomena on spin$^c$ and pin$^-$
manifolds. Let $(\xi, \nabla^{\xi})$ be a real oriented Euclidean
plane bundle, or equivalently a complex line bundle on $M$ with
Euler form $c=e(\xi, \nabla^{\xi})$. For each $8k+4$ dimensional
smooth Riemannian manifold $M$, they obtain that \be \left\{\frac
{\widehat{L}(TM, \nabla^{TM})}{\cosh^2 {({c \over 2})}} \right
\}^{(8k+4)}=8\sum_{r=0}^{k}2^{6k-6r}h_r(T_\CC M, \xi_\CC).\ee In
(1.3), each $h_r(T_\CC M, \xi_\CC)$ is a differential form \h
\left \{\widehat{A}(TM, \nabla^{TM})\mathrm{ch}\left(b_r(T_\CC M,
\xi_\CC), \nabla^{b_r(T_\CC M, \xi_\CC)}\right)\cosh\left({c\over
2}\right)\right\}^{(8k+4)},\e where $\xi_\CC$ is the
complexification of $\xi$ and $b_r(T_\CC M, \xi_\CC)\in
KO(M)\otimes \mathbf{C}, \ 0\leq r\leq k$, can be derived
canonically from $TM$ and $\xi$. (1.3) is a generalization of
Liu's miraculous cancellation formula (1.2) in the sense that when
$c=0$ it exactly gives Liu's formula. Especially on dimension 12,
it gives a generalization of the Alvarez-Gaum\'e-Witten miraculous
cancellation formula,
$$\left\{ {\widehat{L}(TM, \nabla^{TM})\over \cosh^2({c\over 2})}\right\}^{(12)}=\left\{
\left[ 8\widehat{A}(TM, \nabla^{TM}){\rm ch} (
T_{\mathbf{C}}M,\nabla^{T_{\mathbf{C} }M})-32\widehat{A}(TM,
\nabla^{TM})\right.\right.$$
$$- 24\left.\left.\widehat{A}(TM, \nabla^{TM})\left(e^c+e^{-c}-2\right)\right]
\cosh\left({c\over 2}\right)\right\}^{(12)}.$$

As an application of the general cancellation formula (1.3), when
$M$ is closed and spin$^c$ and $B$ is an $8k+2$ dimensional
oriented submanifold of $M$ such that $[B] \in H_{8k+2}(M, Z)$ is
dual to $w_2(TM)$, the authors obtain that
$(\mathrm{sign}(M)-\mathrm{sign}(B\bullet B))$ is divisible by 8
by using the Atiyah-Singer index theorem for spin$^c$ manifolds,
where $(B\bullet B)$ denotes the self-intersection of $B$ in $M$.
Moreover, they [9] show that
\be \begin{split} &\frac{\mathrm{Sig}(M)-\mathrm{Sig}(B\bullet B)}{8}\\
\equiv &\int_M \widehat{A}(TM, \nabla^{TM}) \mathrm{ch}(b_k(T_\CC M+\CC^2-\xi_\CC, \CC^2)) \cosh{\left({c\over 2}\right)}\\
\equiv & \ \mathrm{ind}_2(b_k(TB+\mathbf{R}^2, \mathbf{R}^2)\ \
\mathrm{mod}\ 2,\end{split} \ee which is the analytic version of
the Ochanine congruence obtained in [18]. Formula (1.3) has
interesting applications to study the Ochanine and the Finashin
congruences (cf. [19], [6]). We refer interested readers to [9]
for details. This shows us again how miraculous cancellation
formulas imply divisibility and congruence results.

Looking at these miraculous cancellation formulas as well as the
divisibilities and congruences induced by them, one naturally asks
if there exist more miraculous cancellation formulas like (1.2)
and (1.3) and consequently exist more divisibilities and
congruences for characteristic numbers. We show in this article
that the answer is positive.

To be more precise, still applying modular invariance of
characteristic forms [16], we obtain some interesting twisted
anomaly cancellation formulas. When these new cancellation
formulas are applied to $8k$ and $8k+4$ dimensional closed spin
manifolds, we find some hidden divisibilities of the
characteristic numbers
$\langle\widehat{L}(TM)\mathrm{ch}\left(T_\mathbf{C}M\right),
[M]\rangle$,
$\langle\widehat{L}(TM)\mathrm{ch}\left(T_\mathbf{C}M\otimes
T_\mathbf{C}M \right), [M]\rangle$ and some of their linear
combinations. The divisibilities of the characteristic number
$\langle\widehat{L}(TM)\mathrm{ch}\left(T_\mathbf{C}M\right),
[M]\rangle$ for $8k$ and $8k+4$ dimensional spin manifolds were
already obtained by Hirzebruch [12] by studying elliptic genera.
Our cancellation formulas supply an interesting approach to prove
the Hirzebruch divisibilities. Moreover we are able to construct
examples to show that the Hirzebruch divisibilities are best
possible. On the other hand the divisibilities of the
characteristic number
$\langle\widehat{L}(TM)\mathrm{ch}\left(T_\mathbf{C}M\otimes
T_\mathbf{C}M\right), [M]\rangle$ for $8k$ and $8k+4$ dimensional
spin manifolds induced by our cancellation formulas look new in
the literature. We are also able to construct examples to show
that these divisibilities are best possible. Note that these
characteristic numbers are the indices of the elliptic operators
$d_s\otimes T_\mathbf{C}M$ and $d_s\otimes T_\mathbf{C}M\otimes
T_\mathbf{C}M$ respectively, where $d_s$ is the signature
operator. The twisted signature operator $d_s\otimes
T_\mathbf{C}M$ is already proved to be rigid by using Witten
rigidity theorem, however nobody has been able to give a direct
proof without using it [17]. Our divisibilities are still
applications of modular invariance. It would also be interesting
to find out a direct proof of them.

When we apply our cancellation formulas with a twisted complex
line bundle to spin$^c$ manifolds, we obtain some congruence
results about the characteristic number
$\langle\widehat{L}(TM)\mathrm{ch}\left(T_\mathbf{C}M\right),
[M]\rangle$, which in dimension $8k+4$ give twisted Rokhlin
congruence formulas.

The rest of the article is organized as follows. We list the
twisted anomaly cancellation formulas in Section 2 and postpone
their proofs to Section 5. In Section 3, we apply our twisted
anomaly cancellation formulas to spin manifolds and obtain
divisibilities for the tangent twisted signature. Then in Section
4, the twisted anomaly cancellation formulas are applied to
spin$^c$ manifolds and particularly induce twisted Rokhlin
congruence formulas for $8k+4$ dimensional spin$^c$ manifolds.

\section{Twisted Anomaly Cancellation Formulas}
In this section, we first
present some basic geometric data and then list the twisted miraculous cancellation
formulas.

Let $M$ be a smooth Riemannian manifold. Let $\nabla^{TM}$ be the
associated Levi-Civita connection and $R^{TM}=\nabla^{TM,2}$ the
curvature of $\nabla^{TM}$. $\nabla^{TM}$ extends canonically to a
Hermitian connection $\nabla^{T_\mathbf{C}M}$ on
$T_\mathbf{C}M=TM\otimes\mathbf{C}$.

Let $ \widehat{A}(TM, \nabla^{TM}), \widehat{L}(TM,\nabla^{TM})$ be
the Hirzebruch characteristic forms defined by \be \widehat{A}(TM,
\nabla^{TM})= {\rm det}^{1/2} \left(\frac{\frac{\sqrt{-1}}{4\pi}
R^{TM}}{\sinh \left(\frac{\sqrt{-1}}{4 \pi}R^{TM}\right)}\right),\ee
\be \widehat{L}(TM, \nabla^{TM})= {\rm det}^{1/2}
\left(\frac{\frac{\sqrt{-1}}{2\pi} R^{TM}}{\tanh
\left(\frac{\sqrt{-1}}{4 \pi}R^{TM}\right)}\right).\ee

Let $E,F$ be two Hermitian vector bundles over $M$ carrying
Hermitian connections $\nabla^E, \nabla^F$ respectively. Let
$R^E=\nabla^{E,\ 2}$ (resp. $R^F=\nabla^{F,\ 2}$) be the curvature
of $\nabla^E$ (resp. $\nabla^F$). If we set the formal difference
$G=E-F$, then $G$ carries an induced Hermitian connection
$\nabla^G$ in an obvious sense. We define the associated Chern
character form as
$$ {\rm ch}(G,\nabla^G)={\rm tr}\left[{\rm exp}\left(\frac{\sqrt{-1}}{2
\pi}R^{E}\right)\right]-{\rm tr}\left[{\rm
exp}\left(\frac{\sqrt{-1}}{2 \pi}R^{F}\right)\right].$$

In the rest of the paper, where there will be no confusion about
the Hermitian connection $\nabla^{E}$ on Hermitian vector bundle
$E$, we will write simply $\mathrm{ch}(E)$ for the associated
Chern character form.

Let $\xi$ be a rank two real oriented Euclidean vector bundle, or
equivalently a complex line bundle, over $M$ carrying a Euclidean
connection $\nabla^{\xi}$.

If $E$ is a complex vector bundle over $M$, set
$\widetilde{E}=E-\mathbf{C}^{\mathrm{rk}(E)}. $

Let $q=e^{2 \pi \sqrt{-1}\tau}$ with $\tau \in \mathbf{H}$, the
upper half complex plane.

We introduce four elements (cf. [16], [9]) in
$K(M)[[q^{1\over2}]]$ which consist of formal power series in
$q^{1\over2}$ with coefficients in the $K$-group of $M$, \be
\Theta_1(T_\CC M)=\bigotimes_{n=1}^\infty S_{q^n}(\widetilde{T_\CC
M}) \otimes \bigotimes_{m=1}^\infty \Lambda_{q^m}(\widetilde{T_\CC
M}),\ee

\be \Theta_2(T_\CC M)=\bigotimes_{n=1}^\infty
S_{q^n}(\widetilde{T_\CC M}) \otimes \bigotimes_{m=1}^\infty
\Lambda_{-q^{m-{1\over2}}}(\widetilde{T_\CC M}),\ee

\be \Theta_1(T_\CC M,\xi_\CC)=\bigotimes_{n=1}^\infty
S_{q^n}(\widetilde{T_\CC M}) \otimes \bigotimes_{m=1}^\infty
\Lambda_{q^m}(\widetilde{T_\CC M}-2\widetilde{\xi_\CC})\otimes
\bigotimes_{r=1}^\infty\Lambda_{q^{r-{1\over
2}}}(\widetilde{\xi_\CC})\otimes\bigotimes_{s=1}^\infty\Lambda_{-q^{s-{1\over
2}}}(\widetilde{\xi_\CC}),\ee

\be \Theta_2(T_\CC M,\xi_\CC)=\bigotimes_{n=1}^\infty
S_{q^n}(\widetilde{T_\CC M}) \otimes \bigotimes_{m=1}^\infty
\Lambda_{-q^{m-{1\over2}}}(\widetilde{T_\CC
M}-2\widetilde{\xi_\CC})\otimes
\bigotimes_{r=1}^\infty\Lambda_{q^{r-{1\over
2}}}(\widetilde{\xi_\CC})\otimes\bigotimes_{s=1}^\infty\Lambda_{q^s}(\widetilde{\xi_\CC}).\ee

Recall that for an indeterminate $t$, \be \Lambda_t(E)=\CC
|M+tE+t^2\Lambda^2(E)+\cdots,\ \ \ S_t(E)=\CC |M+tE+t^2
S^2(E)+\cdots, \ee are respectively the total exterior and
symmetric powers of $E$. The following relations hold between
these two operations (cf.
 [2]),
\be S_t(E)=\frac{1}{\Lambda_{-t}(E)},\ \ \
 \Lambda_t(E-F)=\frac{\Lambda_t(E)}{\Lambda_t(F)}.\ee

We can formally expand these four elements into Fourier series in
$q^{1\over2}$,

\be \Theta_1(T_\CC M)=A_0(T_\CC M)+A_1(T_\CC
M)q^{1\over2}+\cdots,\ee

\be \Theta_2(T_\CC M)=B_0(T_\CC M)+B_1(T_\CC M)q^{1\over2}+\cdots,
\ee

\be \Theta_1(T_\CC M,\xi_\CC)=A_0(T_\CC M,\xi_\CC)+A_1(T_\CC
M,\xi_\CC)q^{1\over2}+\cdots, \ee

\be\Theta_2(T_\CC M,\xi_\CC)=B_0(T_\CC M,\xi_\CC)+B_1(T_\CC
M,\xi_\CC)q^{1\over2}+\cdots, \ee where the $A_i$'s and $B_i$'s,
are elements in the semi-group formally generated by Hermitian
vector bundles over $M$. Moreover, they carry canonically induced
Hermitian connections.

Let $c=e(\xi, \nabla^{\xi})$ be the Euler form of $\xi$
canonically associated to $\nabla^{\xi}$.

Now we can state our twisted cancellation formulas and discuss
their applications in the following subsections.

\begin{theorem} For $8k+4$ dimensional smooth Riemannian manifold $M$, the following identity
holds,

\be \left\{\widehat{L}(TM, \nabla^{TM})\mathrm{ch}(T_\CC
M)-16\widehat{L}(TM, \nabla^{TM})\right\}^{(8k+4)}\ee
$$ =2^{14}[\sum_{r=0}^{k-1}(k-r)2^{6(k-r-1)}h_r(T_\CC M)], $$ where
each $h_r(T_\CC M)=\left\{\widehat{A}(TM,
\nabla^{TM})\mathrm{ch}(b_r(T_\CC M)) \right\}^{(8k+4)}, 0\leq r
\leq k,$ and each $b_r(T_\CC M)$ is a canonical integral linear
combination of $B_j(T_\CC M),
 0\leq j \leq r.$ The right hand side is understood as $0$ when
 $k < 1$.
\end{theorem}
In Theorem 2.1, putting $k=0,1,2$, by computing $h_0(T_\CC M)$ and
$h_1(T_\CC M)$ (see (5.18) and (5.19)), we have

\begin{corollary} If $M$ is 4 dimensional, one has
\be \left\{\widehat{L}(TM, \nabla^{TM})\mathrm{ch}(T_\CC
M)-16\,\widehat{L}(TM, \nabla^{TM})\right\}^{(4)}=0.\ee
\end{corollary}

\begin{corollary}If $M$ is 12 dimensional, one has
\be \left\{\widehat{L}(TM, \nabla^{TM})\mathrm{ch}(T_\CC
M)-16\,\widehat{L}(TM,
\nabla^{TM})\right\}^{(12)}=-2^{14}\left\{\widehat{A}(TM,
\nabla^{TM})\right\}^{(12)}.\ee
\end{corollary}

\begin{corollary}If $M$ is 20 dimensional, one has
\be \left\{\widehat{L}(TM, \nabla^{TM})\mathrm{ch}(T_\CC
M)-16\,\widehat{L}(TM, \nabla^{TM})\right\}^{(20)}\ee
$$=2^{14}\left\{\widehat{A}(TM, \nabla^{TM})\mathrm{ch}(T_\CC M)-28\,\widehat{A}(TM, \nabla^{TM})\right\}^{(20)}.$$
\end{corollary}

We also have,

\begin{theorem} For $8k+4$ dimensional smooth Riemannian manifold $M$, the following identity
holds,

\be \left\{\widehat{L}(TM, \nabla^{TM})\mathrm{ch}(T_\CC M\otimes
T_\CC M)-55\widehat{L}(TM, \nabla^{TM})\mathrm{ch}(T_\CC
M)+768\widehat{L}(TM, \nabla^{TM})\right\}^{(8k+4)}\ee
$$ =2^{25}[\sum_{r=0}^{k-2}(k-r)(k-r-1)2^{6(k-r-2)}h_r(T_\CC M)], $$ where
each $h_r(T_\CC M)=\left\{\widehat{A}(TM,
\nabla^{TM})\mathrm{ch}(b_r(T_\CC M)) \right\}^{(8k+4)}, 0\leq r
\leq k,$ and each $b_r(T_\CC M)$ is a canonical integral linear
combination of $B_j(T_\CC M),
 0\leq j \leq r.$ The right hand side is understood as $0$ when
 $k < 2$.
\end{theorem}
In Theorem 2.2, putting $k=0,1,2,3$, by computing $h_0(T_\CC M)$
and $h_1(T_\CC M)$ (see (5.18) and (5.19)), we have

\begin{corollary} If $M$ is 4 dimensional, one has
\be \left\{\widehat{L}(TM, \nabla^{TM})\mathrm{ch}(T_\CC M\otimes
T_\CC M)-55\,\widehat{L}(TM, \nabla^{TM})\mathrm{ch}(T_\CC
M)+768\,\widehat{L}(TM, \nabla^{TM})\right\}^{(4)}=0.\ee
Equivalently, in view of Corollary 2.1, one has \be
\left\{\widehat{L}(TM, \nabla^{TM})\mathrm{ch}(T_\CC M\otimes
T_\CC M)-112\,\widehat{L}(TM, \nabla^{TM})\right\}^{(4)}=0.\ee
\end{corollary}

\begin{corollary} If $M$ is 12 dimensional, one has
\be \left\{\widehat{L}(TM, \nabla^{TM})\mathrm{ch}(T_\CC M\otimes
T_\CC M)-55\,\widehat{L}(TM, \nabla^{TM})\mathrm{ch}(T_\CC
M)+768\,\widehat{L}(TM, \nabla^{TM})\right\}^{(12)}=0.\ee
\end{corollary}

\begin{corollary} If $M$ is 20 dimensional, one has
\be \begin{split} &\left\{\widehat{L}(TM,
\nabla^{TM})\mathrm{ch}(T_\CC M\otimes T_\CC
M)-55\,\widehat{L}(TM, \nabla^{TM})\mathrm{ch}(T_\CC
M)+768\,\widehat{L}(TM,
\nabla^{TM})\right\}^{(20)}\\
=&-2^{26}\left\{\widehat{A}(TM,
\nabla^{TM})\right\}^{(20)}.\end{split}\ee
\end{corollary}

\begin{corollary} If $M$ is 28 dimensional, one has
\be \begin{split} &\left\{\widehat{L}(TM,
\nabla^{TM})\mathrm{ch}(T_\CC M\otimes T_\CC
M)-55\,\widehat{L}(TM, \nabla^{TM})\mathrm{ch}(T_\CC
M)+768\,\widehat{L}(TM,
\nabla^{TM})\right\}^{(28)}\\
=&2^{26}\left\{\widehat{A}(TM, \nabla^{TM})\mathrm{ch}(T_\CC
M)-52\,\widehat{A}(TM, \nabla^{TM})\right\}^{(28)}.
\end{split} \ee
\end{corollary}

For $8k$ dimensional manifolds, we have the following twisted
cancellation formulas.

\begin{theorem} For $8k$ dimensional smooth Riemannian manifold $M$, the following identity
holds, \be \left\{\widehat{L}(TM, \nabla^{TM})\mathrm{ch}(T_\CC
M)\right\}^{(8k)}
=2^{11}[\sum_{r=0}^{k-1}(k-r)2^{6(k-r-1)}z_r(T_\CC M)], \ee where
each $z_r(T_\CC M)=\left\{\widehat{A}(TM,
\nabla^{TM})\mathrm{ch}(d_r(T_\CC M)) \right\}^{(8k)}, 0\leq r
\leq k,$ and each $d_r(T_\CC M)$ is a canonical integral linear
combination of $B_j(T_\CC M), 0\leq j \leq r.$

\end{theorem}

In Theorem 2.3, putting $k=1,2$, by computing $z_0(T_\CC M)$ and
$z_1(T_\CC M)$ (see (5.32) and (5.33)), we have
\begin{corollary} If $M$ is 8 dimensional, one has
\be \left\{\widehat{L}(TM, \nabla^{T_\CC M})\mathrm{ch}(T_\CC
M)\right\}^{(8)}=2048\left\{\widehat{A}(TM,
\nabla^{TM})\right\}^{(8)}.\ee
\end{corollary}

\begin{corollary}If $M$ is 16 dimensional, one has
\be \left\{\widehat{L}(TM, \nabla^{TM})\mathrm{ch}(T_\CC
M)\right\}^{(16)}\ee $$=-2048\left\{\widehat{A}(TM,
\nabla^{TM})\mathrm{ch}(T_\CC M)-48\,\widehat{A}(TM,
\nabla^{TM})\right\}^{(16)}.$$
\end{corollary}

\begin{theorem} For $8k$ dimensional smooth Riemannian manifold $M$, the following identity
holds, \be \begin{split} &\left\{\widehat{L}(TM,
\nabla^{TM})\mathrm{ch}(T_\CC M\otimes T_\CC
M)-23\,\widehat{L}(TM,
\nabla^{TM})\mathrm{ch}(T_\CC M)\right\}^{(8k)}\\
=&2^{22}[\sum_{r=0}^{k-2}(k-r)(k-r-1)2^{6(k-r-2)}z_r(T_\CC M)],
\end{split} \ee where each $z_r(T_\CC M)=\left\{\widehat{A}(TM,
\nabla^{TM})\mathrm{ch}(d_r(T_\CC M)) \right\}^{(8k)}, 0\leq r
\leq k,$ and each $d_r(T_\CC M)$ is a canonical integral linear
combination of $B_j(T_\CC M), 0\leq j \leq r.$ The right hand side
is understood as $0$ when
 $k < 2$.

\end{theorem}

In Theorem 2.4, putting $k=1,2,3$, by computing $z_0(T_\CC M)$ and
$z_1(T_\CC M)$ (see (5.32) and (5.33)), we have
\begin{corollary} If $M$ is 8 dimensional, one has
\be \left\{\widehat{L}(TM, \nabla^{TM})\mathrm{ch}(T_\CC M\otimes
T_\CC M)-23\,\widehat{L}(TM, \nabla^{TM})\mathrm{ch}(T_\CC
M)\right\}^{(8)}=0.\ee Equivalently, in view of Corollary 2.8, one
has \be \left\{\widehat{L}(TM, \nabla^{TM})\mathrm{ch}(T_\CC
M\otimes T_\CC M)\right\}^{(8)}=\left\{23\cdot 2048
\,\widehat{A}(TM, \nabla^{TM})\right\}^{(8)}.\ee
\end{corollary}

\begin{corollary} If $M$ is 16 dimensional, one has
\be \begin{split} &\left\{\widehat{L}(TM,
\nabla^{TM})\mathrm{ch}(T_\CC M\otimes T_\CC
M)-23\,\widehat{L}(TM, \nabla^{TM})\mathrm{ch}(T_\CC
M)\right\}^{(16)}\\
=&2^{23}\left\{\widehat{A}(TM, \nabla^{TM})\right\}^{(16)}.
\end{split}\ee
\end{corollary}

\begin{corollary} If $M$ is 24 dimensional, one has
\be \begin{split} &\left\{\widehat{L}(TM,
\nabla^{TM})\mathrm{ch}(T_\CC M\otimes T_\CC
M)-23\,\widehat{L}(TM, \nabla^{TM})\mathrm{ch}(T_\CC
M)\right\}^{(24)}\\
=&-2^{23}\left\{\widehat{A}(TM, \nabla^{TM})\mathrm{ch}(T_\CC
M)-72\,\widehat{A}(TM, \nabla^{TM})\right\}^{(24)}.
\end{split}\ee
\end{corollary}

If we include the extra complex line bundle $\xi$ into our picture,
we have the following.

\begin{theorem} For $8k+4$ dimensional smooth Riemannian manifold $M$, the following
identity holds, \be \left\{\frac{\widehat{L}(TM,
\nabla^{TM})\left[\mathrm{ch}(T_\CC M, \nabla^{T_\CC
M})-\sinh^2{\left({c\over2}\right)}\mathrm{ch}\left(2\xi_\CC
\oplus\CC^8\right)-16\right]}{\cosh^2{\left({c\over2}\right)}}\right\}^{(8k+4)}\ee
$$=2^{14}[\sum_{r=0}^{k-1}(k-r)2^{6(k-r-1)}h_r(T_\CC M, \xi_\CC)], $$ where each $h_r(T_\CC M, \xi_\CC)=\left\{\widehat{A}(TM,
\nabla^{TM})\mathrm{ch}(b_r(T_\CC M, \xi_\CC))
\cosh{\left({c\over2}\right)}\right\}^{(8k+4)},
 0\leq r \leq k,$ and each $b_r(T_\CC M, \xi_\CC)$ is a canonical integral
linear combination of $B_j(T_\CC M, \xi_\CC),
 0\leq j \leq r.$ The right hand is understood as $0$ when $k <
 1$.

\end{theorem}

For $8k$ dimensional case, we have the following.

\begin{theorem} For $8k$ dimensional smooth Riemannian manifold $M$, the following
identity holds, \be \left\{\frac{\widehat{L}(TM,
\nabla^{TM})\left[\mathrm{ch}(T_\CC M, \nabla^{T_\CC
M})-\sinh^2{\left({c\over2}\right)}\mathrm{ch}\left(2\xi_\CC
\oplus\CC^8\right)\right]}{\cosh^2{\left({c\over2}\right)}}\right\}^{(8k)}\ee
$$=2^{11}[\sum_{r=0}^{k-1}(k-r)2^{6(k-r-1)}z_r(T_\CC M, \xi_\CC)], $$ where each $z_r(T_\CC M, \xi_\CC)=\left\{\widehat{A}(TM,
\nabla^{TM})\mathrm{ch}(d_r(T_\CC M, \xi_\CC))
\cosh{\left({c\over2}\right)}\right\}^{(8k)},
 0\leq r \leq k,$ and each $d_r(T_\CC M, \xi_\CC)$ is a canonical integral
linear combination of $B_j(T_\CC M, \xi_\CC),
 0\leq j \leq r.$
\end{theorem}

\section{Spin Manifolds and Divisibilities of Twisted Signatures} Let $M$ be a closed oriented differential
manifold and $V$ be a complex vector bundle over $M$. Let
$\mathrm{Sig}(M, V)\triangleq \mathrm{Ind}(d_s\otimes V)$ denote
the twisted signature [7]. Let $\mathrm{Sig}(M, T)$ and
$\mathrm{Sig}(M, T \otimes T)$ denote $\mathrm{Sig}(M, T_\CC M)$
and $\mathrm{Sig}(M, T_\CC M \otimes T_\CC M)$ respectively. In
this section, we apply Theorem 2.1, 2.2 and Theorem 2.3, 2.4 to
$8k+4$ and $8k$ dimensional closed spin manifolds respectively to
obtain divisibility results for the twisted signatures
$\mathrm{Sig}(M, T)$ and $\mathrm{Sig}(M, T \otimes T)$. We also
show that our divisibilities are best possible.

\subsection{\bf {$8k+4$ dimensional case}}
According to the generalized Hirzebruch signature formula [7, 11],
when $M$ is an $8k+4$ dimensional closed spin manifold,
integrating both sides of (2.13) against the fundamental class
$[M]$, we have \be \mathrm{Sig}(M,
T)-16\,\mathrm{Sig}(M)=2^{14}\langle
[\sum_{r=0}^{k-1}(k-r)2^{6(k-r-1)}h_r(T_\CC M)], [M]\rangle.\ee
According to the Atiyah-Hirzebruch divisibility [3], we have
\begin{theorem}If $M$ is an $8k+4$ dimensional closed spin
manifold, then  $(\mathrm{Sig}(M, T)-16\,\mathrm{Sig}(M))$ is
divisible by $2^{15}$.
\end{theorem}

Then according to Theorem 3.1 and the Ochanine divisibility [19]
that the signature of $8k+4$ dimensional closed spin manifolds is
divisible by $16$, we see that our twisted anomaly cancellation
formula (2.13) actually implies the Hirzebruch divisibility:

\begin{theorem} [Hirzebruch, {[12]}] If $M$ is an $8k+4$ dimensional closed spin
manifold, then the twisted signature $\mathrm{Sig}(M, T)$ is
divisible by 256.
\end{theorem}

Moreover, we are able to show that the Hirzebruch divisibility is
best possible.
\begin{proposition} 256 is the best possible divisibility of the twisted signature $\mathrm{Sig}(M, T)$ for $8k+4$
dimensional spin manifolds. \end{proposition}

To prove Proposition 3.1, we need the following lemmas.

\begin{lemma} Let $M_1$ and $M_2$ be two closed oriented smooth
manifolds, then one has \be  \mathrm{Sig}(M_1\times M_2,
T)=\mathrm{Sig}(M_1)\mathrm{Sig}(M_2,
T)+\mathrm{Sig}(M_2)\mathrm{Sig}(M_1, T),\ee and \be
\begin{split} &\mathrm{Sig}(M_1\times M_2, T\otimes
T)\\
=&\mathrm{Sig}(M_1)\mathrm{Sig}(M_2, T\otimes
T)+2\mathrm{Sig}(M_1, T)\mathrm{Sig}(M_2,
T)+\mathrm{Sig}(M_2)\mathrm{Sig}(M_1, T\otimes T).\end{split}\ee
\end{lemma}

\begin{proof} Let $p_1: M_1\times M_2\rightarrow M_1,  p_2: M_1\times M_2\rightarrow M_2$ be the two projections. It's not hard to see that \begin{equation*} \begin{split}&\int_{M_1\times M_2}\widehat{L}(M_1\times M_2)\mathrm{ch}(T_C(M_1\times M_2))\\
=&\int_{M_1\times M_2}\widehat{L}(p_1^*(TM_1))\widehat{L}(p_2^*(TM_2))\left(\mathrm{ch}(p_1^*(T_\CC M_1))+\mathrm{ch}(p_2^*(T_\CC M_2))\right)\\
=&\int_{M_1\times M_2}p_1^*\left(\widehat{L}(M_1)\right) p_2^*\left(\widehat{L}(M_2)\mathrm{ch}(T_\CC M_2)\right)+
p_2^*\left(\widehat{L}(M_2)\right) p_1^*\left(\widehat{L}(M_1)\mathrm{ch}(T_\CC M_1)\right)\\
=&\int_{M_1}\widehat{L}(M_1)\int_{M_2}\widehat{L}(M_2)\mathrm{ch}(T_\CC M_2)+\int_{M_2}\widehat{L}(M_2)\int_{M_1}\widehat{L}(M_1)\mathrm{ch}(T_\CC M_1).
\end{split}
\end{equation*}
Thus we have \begin{equation*} \mathrm{Sig}(M_1\times M_2,
T)=\mathrm{Sig}(M_1)\mathrm{Sig}(M_2,
T)+\mathrm{Sig}(M_2)\mathrm{Sig}(M_1, T). \end{equation*}

By similar computations, it's not hard to prove (3.3).
\end{proof}

\begin{lemma} Let $\mathbf{H}P^2$ be the quaternionic projective plane. Then for the $8n$ dimensional
manifold $(\mathbf{H}P^2)^n$, the n-fold product of
$\mathbf{H}P^2$, one has $\mathrm{Sig}((\mathbf{H}P^2)^n)=1$,
$\mathrm{Sig}((\mathbf{H}P^2)^n, T)=0$ and
$\mathrm{Sig}((\mathbf{H}P^2)^n, T\otimes T)=0$, where $n$ is a
positive integer.
\end{lemma}
\begin{proof} Let $\mathbf{H}P^n$ be the quaternionic projective space and $u$ be the
generator of $H^4(\mathbf{H}P^n, \mathbf{Z})$. A theorem of Hirzebruch ([10]) says the total Pontrjagin class of $T\mathbf{H}P^n$
is the following
\be p(\mathbf{H}P^n)=(1+u)^{2n+2}(1+4u)^{-1}=(1+u)^{2n+2}(1-4u+16u^2+\cdots).\ee
In particular, for $\mathbf{H}P^2$, we have $p_1(\mathbf{H}P^2)=2u, p_2(\mathbf{H}P^2)=7u^2.$

By direct computations, Sig$(\mathbf{H}P^2)=\langle -{1\over
45}p_1^2+{7\over 45}p_2, [\mathbf{H}P^2]\rangle=1$. Thus by the
multiplicity of the signature, we have Sig$((\mathbf{H}P^2)^n)=1$,
for $n\in \mathbf{Z}^{+}$.

Also by direct computations, we have
$\widehat{A}(\mathbf{H}P^2)=\langle {1\over 16}({7\over
360}p_1^2-{1\over 90}p_2), [\mathbf{H}P^2]\rangle=0$. Thus by
Corollary 2.8, we have $\mathrm{Sig}(\mathbf{H}P^2, T)=0.$ Keeping
applying (3.2), we obtain that $\mathrm{Sig}((\mathbf{H}P^2)^n,
T)=0,$ for $n\in \mathbf{Z}^{+}$.

By Corollary 2.10, $\mathrm{Sig}(\mathbf{H}P^2, T\otimes
T)=23\,\mathrm{Sig}(\mathbf{H}P^2, T)=0.$ Keeping applying (3.3),
we obtain that $\mathrm{Sig}((\mathbf{H}P^2)^n, T\otimes T)=0,$
for $n\in \mathbf{Z}^{+}$.

\end{proof}

Now we can prove Proposition 3.1 as follows.
\begin{proof}
Let $K$ be a K3-surface. It is well known that Sig$(K)$ is $-16$.
Then by Corollary 2.1, Sig$(K, T)=16\,$Sig$(K)=-256$.

Applying (3.2) and Lemma 3.2 to the $8k+4$ dimensional spin
manifold $K \times(\mathbf{H}P^2)^{k}$, we obtain that
\begin{equation*} \begin{split} &\mathrm{Sig}\left(K\times (\mathbf{H}P^2)^{k}, T\right)\\
=&\mathrm{Sig}(K)\mathrm{Sig}((\mathbf{H}P^2)^{k}, T)+\mathrm{Sig}((\mathbf{H}P^2)^{k})\mathrm{Sig}(K, T)\\
=&\mathrm{Sig}(K, T)\\
=& -256.\end{split}
\end{equation*} This proves Proposition 3.1. \end{proof}

Our twisted anomaly cancellation formula (2.17) implies the
divisibility of the twisted signature $\mathrm{Sig}(M, T\otimes
T)$. According to the generalized Hirzebruch signature formula [7,
11], when $M$ is an $8k+4$ dimensional closed spin manifold,
integrating both sides of (2.17) against the fundamental class
$[M]$, we have \be \begin{split} &\mathrm{Sig}(M, T\otimes
T)-55\,\mathrm{Sig}(M, T)+768\,\mathrm{Sig}(M)\\
=&2^{25}\langle [\sum_{r=0}^{k-2}(k-r)(k-r-1)2^{6(k-r-2)}h_r(T_\CC
M)], [M]\rangle.\end{split} \ee In particular, by (2.19), in
dimension 4, we have \be \mathrm{Sig}(M, T\otimes
T)-112\,\mathrm{Sig}(M)=0.\ee

By the Atiyah-Hirzebruch divisibility [3], (3.5) shows
$\mathrm{Sig}(M, T\otimes T)-55\,\mathrm{Sig}(M,
T)+768\,\mathrm{Sig}(M)$ is divisible by $2^{26}$. Therefore,
according to Theorem 3.2 and the Ochanine divisibility [19], we
obtain

\begin{theorem}Let $M$ be an $8k+4$ dimensional closed spin
manifold, when $\mathrm{dim}M=4$, the twisted signature
$\mathrm{Sig}(M, T\otimes T)$ is divisible by $256\cdot7$; when
$\mathrm{dim}M=8k+4, k \geq 1,$ the the twisted signature
$\mathrm{Sig}(M, T\otimes T)$ is divisible by $256$.
\end{theorem}

Moreover we are also able to show that these divisibilities of
$\mathrm{Sig}(M, T\otimes T)$ are best possible.

\begin{proposition} $256\cdot7$ is the best possible divisibility of the twisted signature $\mathrm{Sig}(M, T\otimes T)$ for $4$
dimensional spin manifolds; $256$ is the best possible
divisibility of the twisted signature $\mathrm{Sig}(M, T\otimes
T)$ for $8k+4$ dimensional spin manifolds, where $k\geq 1$.
\end{proposition}
\begin{proof} Let $K$ be a K3-surface. By (3.6), one has $\mathrm{Sig}(K, T\otimes
T)=112\,\mathrm{Sig}(K)=-256\cdot 7$. This shows $256\cdot7$ is
the best possible divisibility of the twisted signature
$\mathrm{Sig}(M, T\otimes T)$ for $4$ dimensional spin manifolds.

Let $B^8$ be such a Bott manifold, which is $8$ dimensional, spin
with $\widehat{A}(B^8)=1$ and $\mathrm{Sig}(B^8)=0$ [15]. By
Corollary 2.8, $\mathrm{Sig}(B^8, T)=2048\,\widehat{A}(B^8)=2048.$
Then by (2.27), $\mathrm{Sig}(B^8, T\otimes
T)=23\,\mathrm{Sig}(B^8, T)=23\cdot 2048$. Therefore, by (3.3),
\begin{equation*} \begin{split} &\mathrm{Sig}\left(K\times B^8, T\otimes T\right)\\
=&\mathrm{Sig}(K)\mathrm{Sig}(B^8, T\otimes T)+2\mathrm{Sig}(K, T)\mathrm{Sig}(B^8, T)+\mathrm{Sig}(B^8)\mathrm{Sig}(K, T\otimes T)\\
=& -16\cdot 23\cdot 2048+2\cdot(-256)\cdot2048\\
=&-55\cdot 2^{15}.\end{split}
\end{equation*}
Applying (3.3) and Lemma 3.2 to the $8k+4$ dimensional, $k\geq 1$,
spin manifold $(\mathbf{H}P^2)^{k-1}\times K\times B^8$, we have
\begin{equation*}
\begin{split} &\mathrm{Sig}\left((\mathbf{H}P^2)^{k-1}\times K\times B^8, T\otimes T\right)\\
=&\mathrm{Sig}(K\times B^8)\mathrm{Sig}((\mathbf{H}P^2)^{k-1},
T\otimes T)+2\mathrm{Sig}(K\times B^8,
T)\mathrm{Sig}((\mathbf{H}P^2)^{k-1}, T)\\
&+\mathrm{Sig}((\mathbf{H}P^2)^{k-1})\mathrm{Sig}(K\times B^8, T\otimes T)\\
=&\mathrm{Sig}(K\times B^8, T\otimes T)\\
=&-55\cdot 2^{15}.
\end{split}
\end{equation*}
Applying (3.3) and Lemma 3.2 to the $8k+4$ dimensional, $k\geq 1$,
spin manifold $(\mathbf{H}P^2)^{k}\times K$, we have
\begin{equation*}
\begin{split} &\mathrm{Sig}\left((\mathbf{H}P^2)^{k}\times K, T\otimes T\right)\\
=&\mathrm{Sig}(K)\mathrm{Sig}((\mathbf{H}P^2)^{k}, T\otimes
T)+2\mathrm{Sig}(K,
T)\mathrm{Sig}((\mathbf{H}P^2)^{k}, T)\\
&+\mathrm{Sig}((\mathbf{H}P^2)^{k})\mathrm{Sig}(K, T\otimes T)\\
=&\mathrm{Sig}(K, T\otimes T)\\
=&-256\cdot 7.
\end{split}
\end{equation*}
Since the maximal common denominator of $-55\cdot 2^{15}$ and
$-256\cdot 7$ is 256, we see that $256$ is the best possible
divisibility of the twisted signature $\mathrm{Sig}(M, T\otimes
T)$ for $8k+4$ dimensional spin manifolds, where $k\geq 1$.
\end{proof}

\subsection{\bf {$8k$ dimensional case}}

According to the generalized Hirzebruch signature formula [7],
when $M$ is an $8k$ dimensional closed spin manifold, integrating
both sides of (2.23) against the fundamental class $[M]$, we have
\be  \mathrm{Sig}(M, T)=2048\langle
[\sum_{r=0}^{k-1}(k-r)2^{6(k-r-1)}z_r(T_\CC M)], [M]\rangle.\ee
Then according to the Atiyah-Singer index theorem, our anomaly
cancellation formula (2.23) actually implies the Hirzebruch
divisibility:
\begin{theorem} [Hirzebruch, {[12]}] If $M$ is an $8k$ dimensional closed spin
manifold, then $\mathrm{Sig}(M, T)$ is divisible by 2048.
\end{theorem}

\begin{remark} This divisibility looks astonishing since on $8k$ dimensional
closed spin manifold, we can say nothing about the divisibility on
the signature while this twisted signature has so high
divisibility.
\end{remark}

Moreover, we are able to show that the Hirzebruch divisibility is
best possible.
\begin{proposition} $2048$ is the best possible divisibility for the
tangent twisted signature of $8k$ dimensional spin manifolds.
\end{proposition}

\begin{proof}

Let $B^8$ be a Bott manifold as in the proof of Proposition 3.2.
Applying Lemma 3.1 and Lemma 3.2 to the $8k$ dimensional spin
manifold $B^8\times (\mathbf{H}P^2)^{k-1}$, one has
\begin{equation*} \begin{split} &\mathrm{Sig}\left(B^8\times (\mathbf{H}P^2)^{k-1}, T\right)\\
=&\mathrm{Sig}(B^8)\mathrm{Sig}((\mathbf{H}P^2)^{k-1}, T)+\mathrm{Sig}((\mathbf{H}P^2)^{k-1})\mathrm{Sig}(B^8, T)\\
=& \mathrm{Sig}(B^8, T)\\
=&2048.\end{split}
\end{equation*} This proves Proposition 3.3.
\end{proof}

Our twisted anomaly cancellation formula (2.26) implies the
divisibility of the twisted signature $\mathrm{Sig}(M, T\otimes
T)$. According to the generalized Hirzebruch signature formula [7,
11], when $M$ is an $8k$ dimensional closed spin manifold,
integrating both sides of (2.26) against the fundamental class
$[M]$, we have \be \begin{split} &\mathrm{Sig}(M, T\otimes
T)-23\,\mathrm{Sig}(M, T)\\
=&2^{22}\langle [\sum_{r=0}^{k-2}(k-r)(k-r-1)2^{6(k-r-2)}z_r(T_\CC
M)], [M]\rangle.\end{split} \ee In particular, by (2.28), in
dimension 8, we have \be \mathrm{Sig}(M, T\otimes T)-23\cdot
2048\,\widehat{A}(M)=0.\ee

By the Atiyah-Singer index theorem, $\mathrm{Sig}(M, T\otimes
T)-23\,\mathrm{Sig}(M, T)$ is divisible by $2^{22}$. Therefore,
according to Theorem 3.4, we obtain

\begin{theorem}Let $M$ be an $8k$ dimensional closed spin
manifold, when $\mathrm{dim}M=8$, the twisted signature
$\mathrm{Sig}(M, T\otimes T)$ is divisible by $2048\cdot23$; when
$\mathrm{dim}M=8k, k \geq 2,$ the the twisted signature
$\mathrm{Sig}(M, T\otimes T)$ is divisible by $2048$.
\end{theorem}

Moreover, we are also able to show that these divisibilities of
$\mathrm{Sig}(M, T\otimes T)$ are best possible.

\begin{proposition} $2048\cdot23$ is the best possible divisibility of the twisted
signature $\mathrm{Sig}(M, T\otimes T)$ for $8$ dimensional spin
manifolds; $2048$ is the best possible divisibility of the twisted
signature $\mathrm{Sig}(M, T\otimes T)$ for $8k$ dimensional spin
manifolds, where $k\geq 2$.
\end{proposition}
\begin{proof} Let $B^8$ be a Bott manifold as in the proof of Proposition 3.2. $\mathrm{Sig}(B^8, T\otimes
T)=2048\cdot23\,\widehat{A}(B^8)=2048\cdot23$ shows that
$2048\cdot23$ is the best possible divisibility of the twisted
signature $\mathrm{Sig}(M, T\otimes T)$ for $8$ dimensional spin
manifolds.

By (3.3),
\begin{equation*} \begin{split} &\mathrm{Sig}\left(B^8\times B^8, T\otimes T\right)\\
=&\mathrm{Sig}(B^8)\mathrm{Sig}(B^8, T\otimes T)+2\mathrm{Sig}(B^8, T)\mathrm{Sig}(B^8, T)+\mathrm{Sig}(B^8)\mathrm{Sig}(B^8, T\otimes T)\\
=&2\mathrm{Sig}(B^8, T)\mathrm{Sig}(B^8, T)\\
=&2^{23}.\end{split}
\end{equation*}
Applying (3.3) and Lemma 3.2 to the $8k$ dimensional, $k\geq 2$,
spin manifold $(\mathbf{H}P^2)^{k-2}\times B^8\times B^8$, we have
\begin{equation*}
\begin{split} &\mathrm{Sig}\left((\mathbf{H}P^2)^{k-2}\times B^8\times B^8, T\otimes T\right)\\
=&\mathrm{Sig}(B^8\times B^8)\mathrm{Sig}((\mathbf{H}P^2)^{k-2},
T\otimes T)+2\mathrm{Sig}(B^8\times B^8,
T)\mathrm{Sig}((\mathbf{H}P^2)^{k-2}, T)\\
&+\mathrm{Sig}((\mathbf{H}P^2)^{k-2})\mathrm{Sig}(B^8\times B^8, T\otimes T)\\
=&\mathrm{Sig}(B^8\times B^8, T\otimes T)\\
=&2^{23}.
\end{split}
\end{equation*}
Applying (3.3) and Lemma 3.2 to the $8k$ dimensional, $k\geq 2$,
spin manifold $(\mathbf{H}P^2)^{k-1}\times B^8$, we have
\begin{equation*}
\begin{split} &\mathrm{Sig}\left((\mathbf{H}P^2)^{k-1}\times B^8, T\otimes T\right)\\
=&\mathrm{Sig}(B^8)\mathrm{Sig}((\mathbf{H}P^2)^{k-1}, T\otimes
T)+2\mathrm{Sig}(B^8,
T)\mathrm{Sig}((\mathbf{H}P^2)^{k-1}, T)\\
&+\mathrm{Sig}((\mathbf{H}P^2)^{k-1})\mathrm{Sig}(B^8, T\otimes T)\\
=&\mathrm{Sig}(B^8, T\otimes T)\\
=&2048\cdot 23.
\end{split}
\end{equation*}
Since the maximal common denominator of $2^{23}$ and $2048\cdot
23$ is 2048, we see that $2048$ is the best possible divisibility
of the twisted signature $\mathrm{Sig}(M, T\otimes T)$ for $8k$
dimensional spin manifolds, where $k\geq 2$.
\end{proof}

\section{Spin$^c$ Manifolds and Twisted Rokhlin Congruences for Characteristic Numbers } In this
section, we apply Theorem 2.5 and Theorem 2.6 to $8k+4$ and $8k$
dimensional spin$^c$ manifolds respectively to obtain some
congruence results. In particular, for $8k+4$ dimensional spin$^c$
manifolds, we establish twisted Rokhlin congruence formulas.

\subsection{\bf {$8k+4$ dimensional case}}

Let $M$ be an $8k+4$ dimensional closed spin$^c$ manifold and $(\xi,
\nabla^{\xi})$ be a real oriented Euclidean plane bundle, or
equivalently a complex line bundle, over $M$ such that
$w_2(TM)\equiv[e(\xi, \nabla^{\xi})]$ in $H^2(M, \mathbf{Z}_2)$. Let
$B$ be an oriented $8k+2$ dimensional submanifold of $M$ such that
$[B]\in H_{8k+2}(M, \mathbf{Z})$ is dual to $[e(\xi,
\nabla^{\xi})]$. Let $B\bullet B $ denote the self-intersection of
$B$ in $M$ and $N$ be the normal bundle to $B\bullet B $ in $M$.
Applying the Poincar\'e duality, we have \be
\begin{split} \int_M
\frac{\widehat{L}(TM)}{\cosh^2{\left({c\over2}\right)}}=&\int_M
\widehat{L}(TM)-\int_M \widehat{L}(TM)
\frac{\sinh^2{\left({c\over2}\right)}}{\cosh^2{\left({c\over2}\right)}}\\
=& \int_M \widehat{L}(TM)-\int_{B\bullet B} \widehat{L}(B\bullet
B),\end{split} \ee

\be \begin{split}&\int_M \frac{\widehat{L}(TM)\mathrm{ch}(T_\CC
M)}{\cosh^2{\left({c\over2}\right)}}\\
=&\int_M \widehat{L}(TM)\mathrm{ch}(T_\CC M)-\int_M
\widehat{L}(TM)\mathrm{ch}(T_\CC
M)\frac{\sinh^2{\left({c\over2}\right)}}{\cosh^2{\left({c\over2}\right)}}\\
=& \int_M \widehat{L}(TM)\mathrm{ch}(T_\CC M)-\int_{B\bullet
B}\widehat{L}(B\bullet B)\mathrm{ch}(T_\CC(B\bullet B)\oplus N_\CC),
\end{split}\ee
and

\be \int_M  \widehat{L}(TM)\mathrm{ch}\left(2\xi_\CC
\oplus\CC^8\right)\frac{\sinh^2{\left({c\over2}\right)}}{\cosh^2{\left({c\over2}\right)}}=\int_{B\bullet
B}\widehat{L}(B\bullet B)\mathrm{ch}\left(N_\CC \oplus\CC^8\right).
\ee

Therefore, by (2.31), (4.1) to (4.3), we have \be\int_M
\widehat{L}(TM)\mathrm{ch}(T_\CC M)-\int_{B\bullet
B}\widehat{L}(B\bullet B)\mathrm{ch}(T_\CC (B\bullet B)\oplus
2N_\CC\oplus \CC^8) \ee
$$-16\left(\int_M \widehat{L}(TM)-\int_{B\bullet B} \widehat{L}(B\bullet
B)\right)=2^{14}\int_M \sum_{r=0}^{k-1}(k-r)2^{6(k-r-1)}h_r(T_\CC
M,\xi_\CC). $$ Thus one has

\be \begin{split} &{1 \over 128}\left\{\int_M
\widehat{L}(TM)\mathrm{ch}(T_\CC M)-\int_{B\bullet
B}\widehat{L}(B\bullet B)\mathrm{ch}(T_\CC (B\bullet B)\oplus
2N_\CC\oplus \CC^8)\right\} \\
=&\frac{\int_M \widehat{L}(TM)-\int_{B\bullet B}
\widehat{L}(B\bullet B)}{8}+2^{7}\int_M
\sum_{r=0}^{k-1}(k-r)2^{6(k-r-1)}h_r(T_\CC M,\xi_\CC),
\end{split} \ee

and

\be \begin{split}&\frac{1}{2^{14}}\left\{ \int_M
\widehat{L}(TM)\mathrm{ch}(T_\CC M)-\int_{B\bullet
B}\widehat{L}(B\bullet B)\mathrm{ch}(T_\CC (B\bullet B)\oplus
2N_\CC\oplus \CC^8)\right.\\
&\left.-16\left(\int_M \widehat{L}(TM)-\int_{B\bullet B}
\widehat{L}(B\bullet B)\right)\right\}\\
=&\int_M \sum_{r=0}^{k-1}(k-r)2^{6(k-r-1)}h_r(T_\CC M,\xi_\CC)\\
=&\int_M \sum_{r=0}^{k-2}(k-r)2^{6(k-r-1)}h_r(T_\CC
M,\xi_\CC)+\int_M h_{k-1}(T_\CC M,\xi_\CC)\\
=&\int_M \sum_{r=0}^{k-2}(k-r)2^{6(k-r-1)}h_r(T_\CC
M,\xi_\CC)\\
&+\int_M \widehat{A}(TM, \nabla^{TM})\mathrm{ch}(b_{k-1}(T_\CC M,
\xi_\CC)) \cosh{\left({c\over2}\right)}.
\end{split} \ee

From [9, Theorem 3.2] (Theorem 4.2 below), Theorem 2.5, (4.5),
(4.6) and the Atiyah-Singer index theorem for spin$^c$ manifolds,
we obtain that
\begin{theorem} If $M$ is an $8k+4$ dimensional closed spin$^c$ manifold, $\xi$ is a complex line
bundle over $M$ such that $c_1(\xi)\equiv w_2(TM)\in H^2 (M,
\mathbf{Z}_2)$ and $B$ is an oriented $8k+2$ dimensional
submanifold of $M$ such that $[B]\in H_{8k+2}(M, \mathbf{Z})$ is
dual to $c_1(\xi)\in H^2(M, \mathbf{Z})$, then
$$\mathrm{Sig}(M, T_\CC M)-\mathrm{Sig} \left( B\bullet B, T_\CC (B\bullet B)\oplus
2N_\CC\oplus \CC^8 \right)
$$
is divisible by 128 and
$$ \mathrm{Sig}(M, T_\CC M)-\mathrm{Sig} \left( B\bullet B, T_\CC (B\bullet B)\oplus
2N_\CC\oplus \CC^8 \right)-16\left(\mathrm{Sig}(M)-\mathrm{Sig}
\left( B\bullet B\right)\right)$$ is divisible by $2^{14}$.
Moreover, one has \be
\begin{split} &{1 \over 128}\left\{\mathrm{Sig}(M, T_\CC M)-\mathrm{Sig} \left( B\bullet B, T_\CC (B\bullet B)\oplus
2N_\CC\oplus \CC^8 \right)\right\} \\
\equiv&\int_M \widehat{A}(TM, \nabla^{TM})\mathrm{ch}(b_k(T_\CC M,
\xi_\CC)) \cosh{\left({c\over2}\right)}
 \ \ \ \ \ \mathrm{mod}\ 64,
\end{split} \ee
and \be \begin{split}&\frac{1}{2^{14}}\left\{ \mathrm{Sig}(M,
T_\CC M)-\mathrm{Sig} \left( B\bullet B, T_\CC (B\bullet B)\oplus
2N_\CC\oplus \CC^8 \right)-16\left(\mathrm{Sig}(M)-\mathrm{Sig}
\left( B\bullet B\right)\right)\right\}\\ \equiv& \int_M
\widehat{A}(TM, \nabla^{TM})\mathrm{ch}(b_{k-1}(T_\CC M, \xi_\CC))
\cosh{\left({c\over2}\right)}\ \ \ \mathrm{mod}\ 128. \end{split}
\ee

\end{theorem}
Note that we also have the following results:

\begin{theorem} $\mathrm{ (Han-Zhang
\ [9,Theorem\  3.2])}$ The following congruence formula holds,

\be
\begin{split}&\frac{\mathrm{Sig}(M)-\mathrm{Sig}
\left( B\bullet B\right)}{8}\\
\equiv& \int_M \widehat{A}(TM, \nabla^{TM})\mathrm{ch}(b_k(T_\CC
M, \xi_\CC))\cosh{\left({c\over2}\right)} \ \ \mathrm{mod}\ 64.
\end{split}\ee

\end{theorem}

\begin{theorem} $\mathrm{(Liu-Zhang \ [18, Theorem\  4.2], Han-Zhang
\ [9, formula \ 3.5, \ 3.6 \ and \ 3.14]}$ The following congruence
formulas hold, \be
\begin{split}&\int_M \widehat{A}(TM,
\nabla^{TM})\mathrm{ch}(b_r(T_\CC M,
\xi_\CC))\cosh{\left({c\over2}\right)}\\
\equiv& \int_M \widehat{A}(TM, \nabla^{TM})\mathrm{ch}(b_r(T_\CC
M+\CC^2-\xi_\CC, \CC^2))\cosh{\left({c\over2}\right)} \ \
\mathrm{mod}\ 2, \end{split}\ee $0 \leq r \leq k,$ and

\be
\begin{split}&\frac{\mathrm{Sig}(M)-\mathrm{Sig}
\left( B\bullet B\right)}{8}\\
\equiv& \int_M \widehat{A}(TM, \nabla^{TM})\mathrm{ch}(b_k(T_\CC
M+\CC^2-\xi_\CC, \CC^2))\cosh{\left({c\over2}\right)} \ \
\mathrm{mod}\ 2. \end{split}\ee
\end{theorem}

Let $E$ be a real vector bundle over $M$ and $i:B\hookrightarrow
M$ denote the canonical embedding of $B$ in $M$. Then we have ,
\begin{theorem}$\mathrm{(Zhang, [20])}$ The following identity
holds, \be \int_M \widehat{A}(TM,
\nabla^{TM})\mathrm{ch}(E\otimes\CC)\cosh{\left({c\over2}\right)}
\equiv \mathrm{ind}_2(i^*E)\ \ \ \mathrm{mod}\ 2,\ee where
$\mathrm{ind}_2(i^*E)$ is the mod 2 index in the sense of Atiyah
and Singer.

\end{theorem}
Note that Theorem 4.4 only holds for $8k+4$ dimensional spin$^c$
manifolds.

Combining Theorem 4.1, Theorem 4.3 and 4.4, we obtain that
\begin{corollary} The following congruence formulas hold,
\be \begin{split} &{1 \over 128}\left\{\mathrm{Sig}(M, T_\CC
M)-\mathrm{Sig} \left( B\bullet B, T_\CC (B\bullet B)\oplus
2N_\CC\oplus \CC^8 \right)\right\} \\
\equiv& \ \mathrm{ind}_2(b_k(TB+\mathbf{R}^2, \mathbf{R}^2))
 \ \ \ \ \ \mathrm{mod}\ 2,
\end{split} \ee
and \be \begin{split}&\frac{1}{2^{14}}\left\{ \mathrm{Sig}(M,
T_\CC M)-\mathrm{Sig} \left( B\bullet B, T_\CC (B\bullet B)\oplus
2N_\CC\oplus \CC^8 \right)-16\left(\mathrm{Sig}(M)-\mathrm{Sig}
\left( B\bullet B\right)\right)\right\}\\ \equiv& \
\mathrm{ind}_2(b_{k-1}(TB+\mathbf{R}^2, \mathbf{R}^2))
 \ \ \ \ \ \mathrm{mod}\ 2. \end{split} \ee

\end{corollary}

We can regard (4.13) and (4.14) as twisted Rokhlin congruence
formulas for $8k+4$ dimensional spin$^c$ manifolds.
\subsection{\bf {$8k$ dimensional case}}

By similar computations as what we did in the above subsection,
for $8k$ dimensional case, from Theorem 2.6, we have

\be \begin{split} &\int_M \widehat{L}(TM)\mathrm{ch}(T_\CC
M)-\int_{B\bullet B}\widehat{L}(B\bullet B)\mathrm{ch}(T_\CC
(B\bullet B)\oplus 2N_\CC\oplus \CC^8) \\
=&2^{11}\int_M \sum_{r=0}^{k-1}(k-r)2^{6(k-r-1)}z_r(T_\CC M,
\xi_\CC)). \end{split} \ee

\begin{theorem} If $M$ is an $8k$ dimensional closed spin$^c$ manifold, $\xi$ is a complex line
bundle over $M$ such that $c_1(\xi)\equiv w_2(TM)\in H^2(M,
\mathbf{Z}_2)$ and $B$ is an oriented $8k-2$ dimensional submanifold
of $M$ such that $[B]\in H_{8k-2}(M,\mathbf{Z})$ is dual to
$c_1(\xi)\in H^2(M, \mathbf{Z})$. Let $B\bullet B$ be the
self-intersection of $B$ in $M$ and $N$ be the normal bundle of
$B\bullet B$ in $M$,  then
$$\mathrm{Sig}(M, T_\CC M)-\mathrm{Sig} \left( B\bullet B, T_\CC (B\bullet B)\oplus
2N_\CC\oplus \CC^8 \right)
$$
is divisible by 2048. Moreover, one has \be \begin{split} &{1
\over 2048}\left\{\mathrm{Sig}(M, T_\CC M)-\mathrm{Sig} \left(
B\bullet B, T_\CC (B\bullet B)\oplus
2N_\CC\oplus \CC^8 \right) \right)\} \\
\equiv&\int_M z_{k-1}(T_\CC M, \xi_\CC)=\int_M \widehat{A}(TM,
\nabla^{TM})\mathrm{ch}(d_{k-1}(T_\CC M, \xi_\CC))
\cosh{\left({c\over2}\right)} \ \ \ \ \ \mathrm{mod}\ 128.
\end{split} \ee

\end{theorem}

\begin{remark} It's pretty interesting to note that in $8k$ dimensional spin$^c$ case although we can
say nothing about the divisibility of
$(\mathrm{Sig}(M)-\mathrm{Sig}(B\bullet B))$, we do have a very
high divisibility for the twisted version
$$\mathrm{Sig}(M, T_\CC M)-\mathrm{Sig} \left( B\bullet B, T_\CC (B\bullet B)\oplus
2N_\CC\oplus \CC^8 \right),
$$which is even much higher than the $8k+4$ dimensional case.
\end{remark}

\section{Proofs of Twisted Anomaly Cancellation Formulas}

We use the modular invariance method developed in [16, 8, 9] to
prove Theorem 2.1 to 2.6 in this section.

We first recall some necessary knowledge on theta-functions and
modular forms. Then in Section 5.1 we prove Theorem 2.1 to 2.4
together and in Section 5.2 we prove Theorem 2.5 and 2.6 together.

Recall that the four Jacobi theta-functions [5] defined by
infinite multiplications are

 \be \theta(v,\tau)=2q^{1/8}\sin(\pi v)\prod_{j=1}^\infty[(1-q^j)(1-e^{2\pi \sqrt{-1}v}q^j)(1-e^{-2\pi
 \sqrt{-1}v}q^j)], \ee
 \be \theta_1(v,\tau)=2q^{1/8}\cos(\pi v)\prod_{j=1}^\infty[(1-q^j)(1+e^{2\pi \sqrt{-1}v}q^j)(1+e^{-2\pi
 \sqrt{-1}v}q^j)], \ee
 \be \theta_2(v,\tau)=\prod_{j=1}^\infty[(1-q^j)(1-e^{2\pi \sqrt{-1}v}q^{j-1/2})(1-e^{-2\pi
 \sqrt{-1}v}q^{j-1/2})], \ee
\be \theta_3(v,\tau)=\prod_{j=1}^\infty[(1-q^j)(1+e^{2\pi
\sqrt{-1}v}q^{j-1/2})(1+e^{-2\pi \sqrt{-1}v}q^{j-1/2})], \ee where
$q=e^{2\pi \sqrt{-1}\tau}, \tau\in \mathbf{H}$.

They are all holomorphic functions for $(v,\tau)\in \mathbf{C \times
H}$, where $\mathbf{C}$ is the complex plane and $\mathbf{H}$ is the
upper half plane.

Let $\theta^{'}(0,\tau)=\frac{\partial}{\partial
v}\theta(v,\tau)|_{v=0}$, then the following Jacobi identity relates
the four theta-functions garcefully.
\newline

 \begin{proposition}$\mathrm{(Jacobi\ identity, [5, Chapter\  3])}$ The
following identity holds, \be \theta^{'}(0,\tau)=\pi
\theta_1(0,\tau) \theta_2(0,\tau)\theta_3(0,\tau). \ee
\end{proposition}

Let $$ SL_2(\mathbf{Z}):= \left\{\left.\left(\begin{array}{cc}
                                      a&b\\
                                      c&d
                                     \end{array}\right)\right|a,b,c,d\in\mathbf{Z},\ ad-bc=1
                                     \right\}
                                     $$
 as usual be the famous modular group. Let
$$S=\left(\begin{array}{cc}
      0&-1\\
      1&0
\end{array}\right), \ \ \  T=\left(\begin{array}{cc}
      1&1\\
      0&1
\end{array}\right)$$
be the two generators of $ SL_2(\mathbf{Z})$. Their actions on
$\mathbf{H}$ are given by
$$ S:\tau\rightarrow-\frac{1}{\tau}, \ \ \ T:\tau\rightarrow\tau+1.$$

Let
$$ \Gamma_0(2)=\left\{\left.\left(\begin{array}{cc}
a&b\\
c&d
\end{array}\right)\in SL_2(\mathbf{Z})\right|c\equiv0\ \ (\rm mod \ \ 2)\right\},$$

$$ \Gamma^0(2)=\left\{\left.\left(\begin{array}{cc}
a&b\\
c&d
\end{array}\right)\in SL_2(\mathbf{Z})\right|b\equiv0\ \ (\rm mod \ \ 2)\right\}$$
be the two modular subgroup of $SL_2(\mathbf{Z})$. It is known
that the generators of $\Gamma_0(2)$ are $T,ST^2ST$, while the
generators of $\Gamma^0(2)$ are $STS,T^2STS$ (cf. [5]).

If we act theta-functions by $S$ and $T$, the following
transformation formulas hold (cf. [5]), \be
\theta(v,\tau+1)=e^{\pi \sqrt{-1}\over 4}\theta(v,\tau),\ \ \
\theta\left(v,-{1}/{\tau}\right)={1\over\sqrt{-1}}\left({\tau\over
\sqrt{-1}}\right)^{1/2} e^{\pi\sqrt{-1}\tau v^2}\theta\left(\tau
v,\tau\right)\ ;\ee \be \theta_1(v,\tau+1)=e^{\pi \sqrt{-1}\over
4}\theta_1(v,\tau),\ \ \
\theta_1\left(v,-{1}/{\tau}\right)=\left({\tau\over
\sqrt{-1}}\right)^{1/2} e^{\pi\sqrt{-1}\tau v^2}\theta_2(\tau
v,\tau)\ ;\ee \be\theta_2(v,\tau+1)=\theta_3(v,\tau),\ \ \
\theta_2\left(v,-{1}/{\tau}\right)=\left({\tau\over
\sqrt{-1}}\right)^{1/2} e^{\pi\sqrt{-1}\tau v^2}\theta_1(\tau
v,\tau)\ ;\ee \be\theta_3(v,\tau+1)=\theta_2(v,\tau),\ \ \
\theta_3\left(v,-{1}/{\tau}\right)=\left({\tau\over
\sqrt{-1}}\right)^{1/2} e^{\pi\sqrt{-1}\tau v^2}\theta_3(\tau
v,\tau)\ .\ee

Let $\Gamma$ be a subgroup of $SL_2(\mathbf{Z}).$
\begin{definition} A modular form over $\Gamma$ is a holomorphic function $f(\tau)$ on $\mathbf{H}\cup
\{\infty\}$ such that for any
 $$ g=\left(\begin{array}{cc}
             a&b\\
             c&d
             \end{array}\right)\in\Gamma\ ,$$
 the following property holds
 $$f(g\tau):=f(\frac{a\tau+b}{c\tau+d})=\chi(g)(c\tau+d)^kf(\tau), $$
 where $\chi:\Gamma\rightarrow\mathbf{C}^*$ is a character of
 $\Gamma$ and $k$ is called the weight of $f$.
 \end{definition}

If $\Gamma$ is a modular subgroup, let
$\mathcal{M}_\mathbf{{R}}(\Gamma)$ denote the ring of modular
forms over $\Gamma$ with real Fourier coefficients. Writing simply
$\theta_j=\theta_j(0,\tau),\ 1\leq j \leq 3,$ we introduce four
explicit modular forms (cf.\ [16]),
$$ \delta_1(\tau)=\frac{1}{8}(\theta_2^4+\theta_3^4), \ \ \ \
\varepsilon_1(\tau)=\frac{1}{16}\theta_2^4 \theta_3^4\ ,$$
$$\delta_2(\tau)=-\frac{1}{8}(\theta_1^4+\theta_3^4), \ \ \ \
\varepsilon_2(\tau)=\frac{1}{16}\theta_1^4 \theta_3^4\ .$$ They have
the following Fourier expansions in $q^{1/2}$:
$$\delta_1(\tau)={1\over 4}+6q+6q^2+\cdots,\ \ \ \ \varepsilon_1(\tau)={1\over
16}-q+7q^2+\cdots\ , $$
$$\delta_2(\tau)=-{1\over 8}-3q^{1/2}-3q+\cdots,\ \ \ \
\varepsilon_2(\tau)=q^{1/2}+8q+\cdots\ .$$ where the
\textquotedblleft $\cdots$" terms are the higher degree terms, all
of which have integral coefficients. They also satisfy the
transformation laws (cf [14], [16]), \be
\delta_2\left(-\frac{1}{\tau}\right)=\tau^2\delta_1(\tau)\ \ \ \ \
, \ \ \ \ \
\varepsilon_2\left(-\frac{1}{\tau}\right)=\tau^4\varepsilon_1(\tau).\ee

\begin{lemma}$\mathrm{([16])}$ One has that $\delta_1(\tau)\ (resp.\ \varepsilon_1(\tau) ) $
is a modular form of weight $2 \ (resp.\ 4)$ over $\Gamma_0(2)$,
while $\delta_2(\tau) \ (resp.\ \varepsilon_2(\tau))$ is a modular
form of weight $2\ (resp.\ 4)$ over $\Gamma^0(2)$, and moreover
$\mathcal{M}_\mathbf{R}(\Gamma^0(2))=\mathbf{R}[\delta_2(\tau),
\varepsilon_2(\tau)]$.
\end{lemma}

\subsection{Proof of Theorem 2.1 to 2.4} Without loss of generality, we
will adopt the Chern roots formalism as in [16], in the
computations of characteristic forms.

Recall that if $\{w_i\}$ are the fromal Chern roots of a Hermitian
vector bundle $E$ carrying a Hermitian $\nabla^{E}$, then one has
the following formula for the Chern character form of the exterior
power of $E$ [11], \be
\mathrm{ch}(\Lambda_t(E))=\prod_i(1+e^{w_i}t).\ee Let's deal with
$8k+4$ dimensional manifolds first.

For $\tau \in\mathbf{H}$ and $q=e^{2\pi \sqrt{-1}\tau}$, set (cf.
[16]) \be
 P_1(\tau)=\left\{\widehat{L}(TM,\nabla^{TM})
\mathrm{ch}\left(\Theta_1(T_\CC M),\nabla^{\Theta_1(T_\CC
M)}\right)\right\}^{(8k+4)},\ee

\be P_2(\tau)=\left\{ \widehat{A}(TM,\nabla^{TM})
\mathrm{ch}\left(\Theta_2(T_\CC M),\nabla^{\Theta_2(T_\CC
M)}\right)\right\}^{(8k+4)},\ee where $\nabla^{\Theta_i(T_\CC M)},
i=1,2,$ are the Hermitian connections with $q^{j/2}$-coefficients
on $\Theta_i(T_\CC M)$ induced from those on the $A_j(T_\CC M)$'s
and $B_j(T_\CC M)$'s.

Let $\{\pm 2\pi \sqrt{-1}x_j\}$ be the formal Chern roots for
$(T_{\mathbf{C}}M,\nabla^{T_{\mathbf{C}}M})$. In terms of the
theta-functions, we get (cf. [16])\be
P_1(\tau)=2^{4k+2}\left\{\prod_{j=1}^{4k+2}x_j\frac{\theta'(0,\tau)}{\theta(x_j,\tau)}
\frac{\theta_{1}(x_j,\tau)}{\theta_{1}(0,\tau)}\right\}^{(8k+4)},
\ee

\be
P_2(\tau)=\left\{\prod_{j=1}^{4k+2}x_j\frac{\theta'(0,\tau)}{\theta(x_j,\tau)}
\frac{\theta_{2}(x_j,\tau)}{\theta_{2}(0,\tau)}\right\}^{(8k+4)}.\ee

Applying the transformation laws (5.6) to (5.9) for
theta-functions, we see that $P_1(\tau)$ is a modular form of
weight $4k+2$ over $\Gamma_0(2)$; while $P_2(\tau)$ is a modular
form of weight $4k+2$ over $\Gamma^0(2)$. Moreover, the following
identity holds, \be P_1(-1/\tau)={(2\tau)}^{4k+2}P_2(\tau).\ee

Observe that at any point $x\in M$, up to the volume form determined
by the metric on $T_xM$, both $P_i(\tau), i=1,\ 2$, can be viewed as
a power series of $q^{1/2}$ with real Fourier coefficients. Thus,
one can apply Lemma 5.1 to $P_2(\tau)$ to get, at $x$, that \be
P_2(\tau)
 =h_0(T_\CC M)(8\delta_2)^{2k+1}+h_1(T_\CC M)(8\delta_2)^{2k-1}\varepsilon_2
+\cdots+h_k(T_\CC M)(8\delta_2)\varepsilon_2^k ,\ee where each
$h_r(T_\CC M)$, $0\leq r\leq k $, is a real multiple of the volume
form at $x$.

We can show that each $h_r(T_\CC M), 0\leq r\leq k,$ can be
expressed through a canonical integral linear combination of
$\left\{\widehat{A}(TM,\nabla^{TM})\mathrm{ch}\left(B_j(T_\mathbf{C}M),
\nabla^{B_j(T_\mathbf{C}M)}\right) \right\}^{(8k+4)}, 0\leq j\leq
r,$ with coefficients not depending on $x\in M$. As in [16], one
can use the induction method to prove this fact easily by
comparing the coefficients of $q^{j/2}, j\geq 0$, between the two
sides of (5.17). For the consideration of the length of this
paper, we do not give details here
 but only write down the explicit expressions for $h_0(T_\CC M)$ and $h_1(T_\CC M)$ as follows.
\be h_0(T_\CC M)
=-\left\{\widehat{A}(TM,\nabla^{TM})\right\}^{(8k+4)}, \ee

\be h_1(T_\CC
M)=\left\{\widehat{A}(TM,\nabla^{TM})\left[24(2k+1)-\mathrm{ch}\left(B_1(T_\mathbf{C}M),
\nabla^{B_1(T_\mathbf{C}M)}\right)\right] \right\}^{(8k+4)}.\ee

 By (5.16) and
(5.17), we have \be
\begin{split}P_1(\tau)=&2^{4k+2}\frac{1}{\tau^{4k+2}}P_2(-1/\tau)\\
=&2^{4k+2}\frac{1}{\tau^{4k+2}}\Big[h_0(T_\CC
M)\big(8\delta_2(-1/\tau)\big)^{2k+1} +h_1(T_\CC
M)\big(8\delta_2(-1/\tau)\big)^{2k-1}\varepsilon_2(-1/\tau)
+\cdots\\
&+h_k(T_\CC M)\big(8\delta_2(-1/\tau)\big)\big(\varepsilon_2(-1/\tau)\big)^k\Big]\\
=&2^{4k+2}\left[h_0(T_\CC M)(8\delta_1)^{2k+1}+h_1(T_\CC
M)(8\delta_1)^{2k-1}\varepsilon_1 +\cdots+h_k(T_\CC
M)(8\delta_1)\varepsilon_1^k\right]. \end{split}\ee

Expanding $\Theta_1(T_\CC M)$ explicitly, by (2.8) we have \be
\begin{split} \Theta_1(T_\CC M)=&\bigotimes_{n=1}^\infty
S_{q^n}(T_\CC M)\Lambda_{-q^n}(\CC^{8k+4})\\
&\otimes \bigotimes_{m=1}^\infty \Lambda_{q^m}(T_\CC
M)S_{-q^m}(\CC^{8k+4})\\
=&(1+(T_\CC M)q+(S^2T_\CC M)q^2+\cdots)(1+(T_\CC M)q^2+\cdots)\\
&(1-\CC^{8k+4}q+(\Lambda^2\CC^{8k+4})q^2+\cdots)(1-\CC^{8k+4}q^2+\cdots)\\
&(1+(T_\CC M)q+(\Lambda^2T_\CC M)q^2+\cdots)(1+(T_\CC M)q^2+\cdots)\\
&(1-\CC^{8k+4}q+(S^2\CC^{8k+4})q^2+\cdots)(1-\CC^{8k+4}q^2+\cdots)\\
=&(1+2(T_\CC M)q+(T_\CC M\otimes T_\CC M+S^2T_\CC M+\Lambda^2T_\CC
M)q^2+\cdots)\\
&(1+2(T_\CC) Mq^2+\cdots)\\
&(1-2\CC^{8k+4}q+(\CC^{8k+4}\otimes \CC^{8k+4}+S^2\CC^{8k+4}
+\Lambda^2\CC^{8k+4})q^2+\cdots)\\
&(1-2\CC^{8k+4}q^2+\cdots)\\
=&(1+2(T_\CC M)q+2(T_\CC M+T_\CC M\otimes T_\CC M)q^2+\cdots)\\
&(1-2\CC^{8k+4}q+2(\CC^{8k+4}\otimes\CC^{8k+4}-\CC^{8k+4})q^2+\cdots)\\
=&1+2(T_\CC M-\CC^{8k+4})q\\
&+2[-(16k+7)T_\CC M+T_\CC M\otimes T_\CC
M+(8k+4)(8k+3)]q^2+\cdots,\\
\end{split} \ee where the ``$\cdots$" are the terms involving
$q^{j}$'s with $j \geq 3$.

Note that \be
\begin{split} &(8\delta_1)^{2k+1-2r}{\varepsilon_1}^{r}\\
=&(2+48q+48q^2\cdots)^{2k+1-2r}({1\over{16}}-q+7q^2\cdots)^r\\
=&2^{2k+1-6r}[1+24(2k+1-2r)q+24(2k+1-2r)(24k-24r+1)q^2\cdots]\\
&[1-16rq+16(8r^2-r)q^2\cdots]\\
=&2^{2k+1-6r}[1+(48k+24-64r)q+\\
&(1152k^2-3072kr+2048r^2+624k-1024r+24)q^2+\cdots].
\end{split}
 \ee

Therefore, by (5.12), (5.21) and (5.22), setting $q=0$ in (5.20),
we get the result of Liu ([16]) \be
\left\{\widehat{L}(TM,\nabla^{TM})\right\}^{(8k+4)}
=8\sum_{r=0}^{k}2^{6k-6r}h_r(T_\CC M).\ee

On the other hand, by (5.12), (5.21) and (5.22), comparing the
coefficients of $q$ in (5.20), we have \be \begin{split} &\left\{
\widehat{L}(TM,\nabla^{TM})\mathrm{ch}\left(2(T_\CC
M)-2(8k+4)\right)
  \right\}^{(8k+4)}\\
=&8\sum_{r=0}^{k}h_r(T_\CC M)2^{6k-6r}(48k+24-64r).\end{split} \ee
Thus \be \begin{split} &\left\{
\widehat{L}(TM,\nabla^{TM})\mathrm{ch}\left((T_\CC
M)-(8k+4)\right)
  \right\}^{(8k+4)}\\
=&4\sum_{r=0}^{k}h_r(T_\CC M)2^{6k-6r}(48k+24-64r)\\
=&12\left(8\sum_{r=0}^{k}2^{6k-6r}h_r(T_\CC
M)\right)-8k\left(8\sum_{r=0}^{k}2^{6k-6r}h_r(T_\CC
M)\right)\\
&+\sum_{r=0}^{k}h_r(T_\CC M)2^{6k-6r}(256k-256r)\\
=&12\left(8\sum_{r=0}^{k}2^{6k-6r}h_r(T_\CC
M)\right)-8k\left(8\sum_{r=0}^{k}2^{6k-6r}h_r(T_\CC
M)\right)\\
&+256\cdot 2^6\sum_{r=0}^{k-1}(k-r)h_{r}(T_\CC
M)2^{6k-6r-6}.\end{split} \ee

Combining (5.23) and (5.25), one has \be
\begin{split}&\left\{\widehat{L}(TM, \nabla^{T_\CC M})\mathrm{ch}(T_\CC M)-16\widehat{L}(TM,
\nabla^{TM})\right\}^{(8k+4)}\\
=&2^{8}[\sum_{r=0}^{k}(k-r)2^{6(k-r)}h_r(T_\CC M)]\\
=&2^{14}[\sum_{r=0}^{k-1}(k-r)2^{6(k-r-1)}h_r(T_\CC M)],
\end{split}\ee which is just (2.13).

Furthermore, by (5.12), (5.21) and (5.22), comparing the
coefficients of $q^2$ in (5.20), we have \begin{equation*}
\begin{split} &\left\{
\widehat{L}(TM,\nabla^{TM})\mathrm{ch}\left(2[-(16k+7)T_\CC
M+T_\CC M\otimes T_\CC M+(8k+4)(8k+3)]\right)
  \right\}^{(8k+4)}\\
=&8\sum_{r=0}^{k}h_r(T_\CC
M)2^{6k-6r}(1152k^2-3072kr+2048r^2+624k-1024r+24).\end{split}
\end{equation*}
Thus, combining (5.23) and (5.26), we have
\begin{equation*}
\begin{split} &\left\{
\widehat{L}(TM,\nabla^{TM})\mathrm{ch}\left([-(16k+7)T_\CC M+T_\CC
M\otimes T_\CC M+(8k+4)(8k+3)]\right)
  \right\}^{(8k+4)}\\
=&32\sum_{r=0}^{k}h_r(T_\CC
M)2^{6k-6r}(144k^2-384kr+256r^2+78k-128r+3)\\
=&-16k\cdot2^{8}[\sum_{r=0}^{k}(k-r)2^{6(k-r)}h_r(T_\CC
M)]+48\cdot2^{8}[\sum_{r=0}^{k}(k-r)2^{6(k-r)}h_r(T_\CC
M)]\\
&+(64k^2-200k+12)\cdot8\sum_{r=0}^{k}2^{6k-6r}h_r(T_\CC
M)\\
&+2^{13}\sum_{r=0}^{k}(k-r)(k-r-1)2^{6(k-r)}h_r(T_\CC
M)\\
=&-16k\left\{\widehat{L}(TM, \nabla^{T_\CC M})\mathrm{ch}(T_\CC
M)-16\widehat{L}(TM, \nabla^{TM})\right\}^{(8k+4)}\\
&+48\left\{\widehat{L}(TM, \nabla^{T_\CC M})\mathrm{ch}(T_\CC
M)-16\widehat{L}(TM, \nabla^{TM})\right\}^{(8k+4)}\\
&+(64k^2-200k+12)\left\{\widehat{L}(TM, \nabla^{T_\CC
M})\right\}^{(8k+4)}\\
&+2^{13}\sum_{r=0}^{k}(k-r)(k-r-1)2^{6(k-r)}h_r(T_\CC
M).\end{split}
\end{equation*}
Therefore by above computations, we have
\begin{equation*}
\begin{split} &\left\{\widehat{L}(TM, \nabla^{TM})\mathrm{ch}(T_\CC M\otimes
T_\CC M)-55\widehat{L}(TM, \nabla^{TM})\mathrm{ch}(T_\CC
M)+768\widehat{L}(TM, \nabla^{TM})\right\}^{(8k+4)}\\
=&2^{13}\sum_{r=0}^{k}(k-r)(k-r-1)2^{6(k-r)}h_r(T_\CC M)\\
=&2^{25}\sum_{r=0}^{k-2}(k-r)(k-r-1)2^{6(k-r-2)}h_r(T_\CC M),
\end{split} \end{equation*}
which is just (2.17).

To prove Theorem 2.2 for $8k$ dimensional case, similarly we set
(cf. [16]) \be
 P_1(\tau)=\left\{\widehat{L}(TM,\nabla^{TM})
\mathrm{ch}\left(\Theta_1(T_\CC M),\nabla^{\Theta_1(T_\CC
M)}\right)\right\}^{(8k)},\ee

\be P_2(\tau)=\left\{ \widehat{A}(TM,\nabla^{TM})
\mathrm{ch}\left(\Theta_2(T_\CC M),\nabla^{\Theta_2(T_\CC
M)}\right)\right\}^{(8k)}.\ee Then one similarly finds that
$P_1(\tau)$ is a modular form of weight $4k$ over $\Gamma_0(2)$;
while $P_2(\tau)$ is a modular form of weight $4k$ over
$\Gamma^0(2)$ and \be P_1(-1/\tau)={(2\tau)}^{4k}P_2(\tau).\ee

This time, applying Lemma 5.1, we have \be P_2(\tau)
 =z_0(T_\CC M)(8\delta_2)^{2k}+z_1(T_\CC M)(8\delta_2)^{2k-2}\varepsilon_2
+\cdots+z_k(T_\CC M)\varepsilon_2^k .\ee And thus

 \be P_1(\tau)
 =2^{4k}\left[z_0(T_\CC M)(8\delta_1)^{2k}+z_1(T_\CC M)(8\delta_1)^{2k-2}\varepsilon_1
+\cdots+z_k(T_\CC M)\varepsilon_1^k \right].\ee

Now we have,

\be z_0(T_\CC M) =\left\{\widehat{A}(TM,\nabla^{TM})\right\}^{(8k)},
\ee

\be z_1(T_\CC
M)=-\left\{\widehat{A}(TM,\nabla^{TM})\left[48k-\mathrm{ch}\left(B_1(T_\mathbf{C}M),
\nabla^{B_1(T_\mathbf{C}M)}\right)\right] \right\}^{(8k)}.\ee

Expanding $\Theta_1(T_\CC M)$ explicitly, we have \be
\begin{split} \Theta_1(T_\CC M)=&\bigotimes_{n=1}^\infty
S_{q^n}(T_\CC M)\Lambda_{-q^n}(\CC^{8k})\\
&\otimes \bigotimes_{m=1}^\infty \Lambda_{q^m}(T_\CC
M)S_{-q^m}(\CC^{8k})\\
=&(1+(T_\CC M)q+(S^2T_\CC M)q^2+\cdots)(1+(T_\CC M)q^2+\cdots)\\
&(1-\CC^{8k}q+(\Lambda^2\CC^{8k})q^2+\cdots)(1-\CC^{8k}q^2+\cdots)\\
&(1+(T_\CC M)q+(\Lambda^2T_\CC M)q^2+\cdots)(1+(T_\CC M)q^2+\cdots)\\
&(1-\CC^{8k}q+(S^2\CC^{8k})q^2+\cdots)(1-\CC^{8k}q^2+\cdots)\\
=&(1+2(T_\CC M)q+(T_\CC M\otimes T_\CC M+S^2T_\CC M+\Lambda^2T_\CC
M)q^2+\cdots)\\
&(1+2(T_\CC) Mq^2+\cdots)\\
&(1-2\CC^{8k}q+(\CC^{8k}\otimes \CC^{8k}+S^2\CC^{8k}
+\Lambda^2\CC^{8k})q^2+\cdots)\\
&(1-2\CC^{8k}q^2+\cdots)\\
=&(1+2(T_\CC M)q+2(T_\CC M+T_\CC M\otimes T_\CC M)q^2+\cdots)\\
&(1-2\CC^{8k}q+2(\CC^{8k}\otimes\CC^{8k}-\CC^{8k})q^2+\cdots)\\
=&1+2(T_\CC M-\CC^{8k})q\\
&+2[-(16k-1)T_\CC M+T_\CC M\otimes T_\CC
M+8k(8k-1)]q^2+\cdots,\\
\end{split} \ee where the ``$\cdots$" are the terms involving
$q^{j}$'s with $j \geq 3$. Note that \be
\begin{split} &(8\delta_1)^{2k-2r}{\varepsilon_1}^{r}\\
=&(2+48q+48q^2\cdots)^{2k-2r}({1\over{16}}-q+7q^2\cdots)^r\\
=&2^{2k-6r}[1+24(2k-2r)q+24(k-r)(48k-48r-22)q^2\cdots]\\
&[1-16rq+16(8r^2-r)q^2\cdots]\\
=&2^{2k-6r}[1+(48k-64r)q+\\
&(1152k^2-3072kr+2048r^2-528k+512r)q^2+\cdots].
\end{split}
 \ee

Therefore by (5.27), (5.34) and (5.35), comparing the constant
terms of both sides of (5.31), we get the result of Liu ([16]) \be
\left\{\widehat{L}(TM,\nabla^{TM})\right\}^{(8k)}
=\sum_{r=0}^{k}2^{6k-6r}z_r(T_\CC M).\ee By (5.27), (5.34) and
(5.35), comparing the coefficients of $q$ of both sides of (5.31),
we have \be \left\{
\widehat{L}(TM,\nabla^{TM})\mathrm{ch}\left(2(T_\CC M,
\nabla^{T_\CC M})-2(8k)\right)
  \right\}^{(8k)}=\sum_{r=0}^{k}2^{6k-6r}(48k-64r)z_r(T_\CC M).
\ee

Thus \be \begin{split} &\left\{
\widehat{L}(TM,\nabla^{TM})\mathrm{ch}\left((T_\CC M, \nabla^{T_\CC
M})-8k\right)
  \right\}^{(8k)}\\
=& \sum_{r=0}^{k}2^{6k-6r}(24k-32r)z_r(T_\CC M)\\
=&\ -8k\left(\sum_{r=0}^{k}2^{6k-6r}z_r(T_\CC M)\right)+
\sum_{r=0}^{k}2^{6k-6r}(32k-32r)z_r(T_\CC M)\\
=&\ -8k\left(\sum_{r=0}^{k}2^{6k-6r}z_r(T_\CC M)\right)+
32\cdot2^6\sum_{r=0}^{k-1}(k-r)2^{6k-6r-6}z_r(T_\CC M).
\end{split} \ee

Combining (5.36) and (5.38), we get \be
\begin{split} &\left\{\widehat{L}(TM, \nabla^{TM})\mathrm{ch}(T_\CC M,
\nabla^{T_\CC
M})\right\}^{(8k)}\\
=&2^{5}[\sum_{r=0}^{k}(k-r)2^{6(k-r)}z_r(T_\CC M)]\\
=&2^{11}[\sum_{r=0}^{k-1}(k-r)2^{6(k-r-1)}z_r(T_\CC M)],
\end{split}\ee which is just (2.23).

By (5.27), (5.34) and (5.35), comparing the coefficients of $q^2$
of both sides of (5.31), we have \begin{equation*}
\begin{split} &\left\{
\widehat{L}(TM,\nabla^{TM})\mathrm{ch}\left(2[-(16k-1)T_\CC
M+T_\CC M\otimes T_\CC M+8k(8k-1)]\right)
  \right\}^{(8k)}\\
=&\sum_{r=0}^{k}z_r(T_\CC
M)2^{6k-6r}(1152k^2-3072kr+2048r^2-528k+512r).\end{split}
\end{equation*}

Thus, combining (5.36) and (5.39), we have
\begin{equation*}
\begin{split} &\left\{
\widehat{L}(TM,\nabla^{TM})\mathrm{ch}\left([-(16k-1)T_\CC M+T_\CC
M\otimes T_\CC M+8k(8k-1)]\right)
  \right\}^{(8k)}\\
=&\sum_{r=0}^{k}z_r(T_\CC
M)2^{6k-6r}(576k^2-1536kr+1024r^2-264k+256r)\\
=&-16k\cdot2^{5}[\sum_{r=0}^{k}(k-r)2^{6(k-r)}z_r(T_\CC
M)]+24\cdot2^{5}[\sum_{r=0}^{k}(k-r)2^{6(k-r)}z_r(T_\CC
M)]\\
&+(64k^2-8k)\sum_{r=0}^{k}2^{6k-6r}z_r(T_\CC
M)\\
&+2^{10}\sum_{r=0}^{k}(k-r)(k-r-1)2^{6(k-r)}z_r(T_\CC
M)\\
=&-16k\left\{\widehat{L}(TM, \nabla^{T_\CC M})\mathrm{ch}(T_\CC
M)\right\}^{(8k)}\\
&+24\left\{\widehat{L}(TM, \nabla^{T_\CC M})\mathrm{ch}(T_\CC
M)\right\}^{(8k)}\\
&+(64k^2-8k)\left\{\widehat{L}(TM, \nabla^{T_\CC
M})\right\}^{(8k+4)}\\
&+2^{10}\sum_{r=0}^{k}(k-r)(k-r-1)2^{6(k-r)}z_r(T_\CC
M).\end{split}
\end{equation*}
Therefore by above computations, we have
\begin{equation*}
\begin{split} &\left\{\widehat{L}(TM, \nabla^{TM})\mathrm{ch}(T_\CC M\otimes
T_\CC M)-23\widehat{L}(TM, \nabla^{TM})\mathrm{ch}(T_\CC
M)\right\}^{(8k)}\\
=&2^{10}\sum_{r=0}^{k}(k-r)(k-r-1)2^{6(k-r)}z_r(T_\CC M)\\
=&2^{22}\sum_{r=0}^{k-2}(k-r)(k-r-1)2^{6(k-r-2)}z_r(T_\CC M),
\end{split} \end{equation*}
which is just (2.26).

\subsection{Proof of Theorem 2.5 and 2.6}

The proof for the cases with the extra complex line bundle $\xi$
involved is similar to the above proof.

For $\tau \in\mathbf{H}$ and $q=e^{2\pi \sqrt{-1}\tau}$, set (cf.
[8, 9]) \be
 P_1(\xi_\CC,
 \tau)=\left\{\frac{\widehat{L}(TM,\nabla^{TM})}{\cosh^2{\left({c\over2}\right)}}
\mathrm{ch}\left(\Theta_1(T_\CC M, \xi_\CC),\nabla^{\Theta_1(T_\CC
M, \xi_\CC)}\right)\right\}^{(8k+4)},\ee

\be P_2(\xi_\CC, \tau)=\left\{ \widehat{A}(TM,\nabla^{TM})
\mathrm{ch}\left(\Theta_2(T_\CC M, \xi_\CC),\nabla^{\Theta_2(T_\CC
M,
\xi_\CC)}\right)\cosh{\left({c\over2}\right)}\right\}^{(8k+4)},\ee
where $\nabla^{\Theta_i(T_\CC M, \xi_\CC)}, i=1,2,$ are the
Hermitian connections with $q^{j/2}$-coefficients on
$\Theta_i(T_\CC M, \xi_\CC)$ induced from those on the $A_j(T_\CC
M, \xi_\CC)$'s and $B_j(T_\CC M, \xi_\CC)$'s.

Let $\{\pm 2\pi \sqrt{-1}x_j\}$ be the formal Chern roots for
$(T_{\mathbf{C}}M,\nabla^{T_{\mathbf{C}}M}),c=2\pi\sqrt{-1}u.$ In
terms of the theta-functions, we get (cf. [9]) \be P_1(\xi_\CC,
\tau)=2^{4k+2}\left\{\left(\prod_{j=1}^{4k+2}x_j\frac{\theta'(0,\tau)}{\theta(x_j,\tau)}
\frac{\theta_{1}(x_j,\tau)}{\theta_{1}(0,\tau)}\right)\frac{\theta_{1}^2(0,\tau)}{\theta_{1}^2(u,\tau)}
\frac{\theta_{3}(u,\tau)}{\theta_{3}(0,\tau)}
\frac{\theta_{2}(u,\tau)}{\theta_{2}(0,\tau)} \right\}^{(8k+4)},
\ee

\be P_2(\xi_\CC,
\tau)=\left\{\left(\prod_{j=1}^{4k+2}x_j\frac{\theta'(0,\tau)}{\theta(x_j,\tau)}
\frac{\theta_{2}(x_j,\tau)}{\theta_{2}(0,\tau)}\right)
\frac{\theta_2^2(0,\tau)}{\theta_2^2(u,\tau)}
\frac{\theta_3(u,\tau)}{\theta_3(0,\tau)}
\frac{\theta_1(u,\tau)}{\theta_1(0,\tau)} \right\}^{(8k+4)}.\ee

Applying the transformation laws (5.6) to (5.9) for
theta-functions, we still see that ([9]) $P_1(\xi_\CC, \tau)$ is a
modular form of weight $4k+2$ over $\Gamma_0(2)$; while
$P_2(\xi_\CC, \tau)$ is a modular form of weight $4k+2$ over
$\Gamma^0(2)$. Moreover, the following identity holds, \be
P_1(\xi_\CC, -1/\tau)={(2\tau)}^{4k+2}P_2(\xi_\CC, \tau).\ee Then
similar to (5.17), we have \be P_2(\xi_\CC, \tau)
 =h_0(T_\CC M, \xi_\CC)(8\delta_2)^{2k+1}+h_1(T_\CC M, \xi_\CC)(8\delta_2)^{2k-1}\varepsilon_2
+\cdots+h_k(T_\CC M, \xi_\CC)(8\delta_2)\varepsilon_2^k ,\ee where
each $h_r(T_\CC M, \xi_\CC), 0\leq r\leq k,$ can be expressed
through a canonical integral linear combination of
$\left\{\widehat{A}(TM,\nabla^{TM})\mathrm{ch}\left(B_j(T_\mathbf{C}M,\xi_\CC),
\nabla^{B_j(T_\mathbf{C}M,\xi_\CC)}\right)\cosh{\left({c\over2}\right)}
\right\}^{(8k+4)}, 0\leq j\leq r.$

Explicitly, one has ([9])

\be h_0(T_\CC M, \xi_\CC)
=-\left\{\widehat{A}(TM,\nabla^{TM})\cosh{\left({c\over2}\right)}\right\}^{(8k+4)},
\ee

\be h_1(T_\CC M,
\xi_\CC)=\left\{\widehat{A}(TM,\nabla^{TM})\left[24(2k+1)-\mathrm{ch}\left(B_1(T_\mathbf{C}M,
\xi_\CC)\right)\right]
\cosh{\left({c\over2}\right)}\right\}^{(8k+4)}.\ee

By (5.44) and (5.45), we have \be \begin{split} P_1(\xi_\CC, \tau)
=&2^{4k+2}\left[h_0(T_\CC M, \xi_\CC)(8\delta_1)^{2k+1}+h_1(T_\CC M,
\xi_\CC)(8\delta_1)^{2k-1}\varepsilon_1 \right.\\
&+\left.\cdots+h_k(T_\CC M,
\xi_\CC)(8\delta_1)\varepsilon_1^k\right]. \end{split} \ee

Let's explicitly expand $\Theta_1(T_\CC M,\xi_\CC)$. By (2.5) and
(2.8), we have \be \begin{split}\Theta_1(T_\CC
M,\xi_\CC)=&\bigotimes_{n=1}^\infty S_{q^n}(\widetilde{T_\CC M})
\otimes \bigotimes_{m=1}^\infty \Lambda_{q^m}(\widetilde{T_\CC
M}-2\widetilde{\xi_\CC})\\ &\otimes
\bigotimes_{r=1}^\infty\left(\Lambda_{q^{r-{1\over
2}}}(\xi_\CC)S_{-q^{r-{1\over
2}}}(\CC^2)\right)\otimes\left(\bigotimes_{s=1}^\infty\Lambda_{-q^{s-{1\over
2}}}(\xi_\CC)S_{q^{s-{1\over 2}}}(\CC^2)\right)\\
=&(1+(T_\CC M-(8k+4))q)\otimes(1+(T_\CC M-(8k+4)-2\xi_\CC+4)q)\\
&\otimes(1+(\xi_\CC) q^{1\over2}+(\xi_\CC\wedge\xi_\CC)q)\otimes(1-\CC^2q^{1\over2}+\CC^3q)\\
&\otimes(1-(\xi_\CC) q^{1\over2}+(\xi_\CC\wedge\xi_\CC)q)\otimes(1+\CC^2q^{1\over2}+\CC^3q)+\cdots\\
=&1+[2(T_\CC M-(8k+4)-\xi_\CC+2)-(\xi_\CC\otimes\xi_\CC-2\xi_\CC\wedge \xi_\CC-\CC^2)]q+\cdots,\\
\end{split}\ee
where the ``$\cdots$" terms are the terms involving
$q^{j\over2}$'s with $j \geq 3.$

By (5.22), (5.40) and (5.49), setting $q=0$ in (5.48), we have
([9])\be
\left\{\frac{\widehat{L}(TM,\nabla^{TM})}{\cosh^2{\left({c\over2}\right)}}\right\}^{(8k+4)}
=8\sum_{r=0}^{k}2^{6k-6r}h_r(T_\CC M, \xi_\CC).\ee

On the other hand, by (5.22), (5.40) and (5.49), comparing the
coefficients of $q$ in (5.48), we have \be
\begin{split}&\left\{
\frac{\widehat{L}(TM,\nabla^{TM})}{\cosh^2{\left({c\over2}\right)}}\mathrm{ch}\left(2(T_\CC
M)-2(8k+4)-2(\xi_\CC-2)-(\xi_\CC\otimes\xi_\CC-2\xi_\CC\wedge
\xi_\CC-\CC^2)\right)
  \right\}^{(8k+4)}\\
  =&8\sum_{r=0}^{k}2^{6k-6r}(48k+24-64r)h_r(T_\CC M, \xi_\CC).\\
  \end{split}
\ee

Thus similar to (5.25), we have

\be
\begin{split}&\left\{
\frac{\widehat{L}(TM,\nabla^{TM})}{\cosh^2{\left({c\over2}\right)}}\mathrm{ch}\left(T_\CC
M-(8k+4)-(\xi_\CC-2)-{1\over2}(\xi_\CC\otimes\xi_\CC-2\xi_\CC\wedge
\xi_\CC-\CC^2)\right)
  \right\}^{(8k+4)}\\
  =&12\left(8\sum_{r=0}^{k}2^{6k-6r}h_r(T_\CC
M, \xi_\CC)\right)-8k\left(8\sum_{r=0}^{k}2^{6k-6r}h_r(T_\CC
M, \xi_\CC)\right)\\
&+256\cdot 2^6\sum_{r=0}^{k-1}(k-r)h_{r}(T_\CC M,
\xi_\CC)2^{6k-6r-6}.\end{split} \ee

Note that \be \begin{split} &
\frac{\widehat{L}(TM,\nabla^{TM})}{\cosh^2{\left({c\over2}\right)}}\mathrm{ch}\left((\xi_\CC-2)+{1\over2}(\xi_\CC\otimes\xi_\CC-2\xi_\CC\wedge
\xi_\CC-\CC^2)\right)\\
=&
\frac{\widehat{L}(TM,\nabla^{TM})}{\cosh^2{\left({c\over2}\right)}}\left((e^{c}+e^{-c}-2)+{1\over2}((e^{c}+e^{-c})^2-4)\right)\\
=&
\frac{\widehat{L}(TM,\nabla^{TM})}{\cosh^2{\left({c\over2}\right)}}\sinh^2{\left({c\over2}\right)}(2(e^{c}+e^{-c})+8)\\
=&
\frac{\widehat{L}(TM,\nabla^{TM})}{\cosh^2{\left({c\over2}\right)}}\sinh^2{\left({c\over2}\right)}\mathrm{ch}(2\xi_\CC\oplus\CC^8).
\end{split}\ee

Therefore combining (5.50),(5.52) and (5.53), one has  \be
\left\{\frac{\widehat{L}(TM, \nabla^{TM})\left[\mathrm{ch}(T_\CC
M)-\sinh^2{({c\over2})}\mathrm{ch}\left(2\xi_\CC\oplus\CC^8\right)-16\right]}{\cosh^2{({c\over2})}}\right\}^{(8k+4)}\ee
$$=2^{14}[\sum_{r=0}^{k-1}(k-r)2^{6(k-r-1)}h_r(T_\CC M, \xi_\CC)], $$ which is just (2.31).

To prove Theorem 2.6 for $8k$ dimensional case, similarly we set
([9]) \be
 P_1(\xi_\CC,
 \tau)=\left\{\frac{\widehat{L}(TM,\nabla^{TM})}{\cosh^2{\left({c\over2}\right)}}
\mathrm{ch}\left(\Theta_1(T_\CC M, \xi_\CC),\nabla^{\Theta_1(T_\CC
M, \xi_\CC)}\right)\right\}^{(8k)},\ee

\be P_2(\xi_\CC, \tau)=\left\{ \widehat{A}(TM,\nabla^{TM})
\mathrm{ch}\left(\Theta_2(T_\CC M, \xi_\CC),\nabla^{\Theta_2(T_\CC
M,
\xi_\CC)}\right)\cosh{\left({c\over2}\right)}\right\}^{(8k)}.\ee

Still playing the same game, we see that $P_1(\xi_\CC, \tau)$ is a
modular form of weight $4k$ over $\Gamma_0(2)$; while
$P_2(\xi_\CC, \tau)$ is a modular form of weight $4k$ over
$\Gamma^0(2)$ and one has the following identities, \be
P_1(-1/\tau)={(2\tau)}^{4k}P_2(\tau),\ee \be P_2(\tau)
 =z_0(T_\CC M, \xi_\CC)(8\delta_2)^{2k}+z_1(T_\CC M, \xi_\CC)(8\delta_2)^{2k-2}\varepsilon_2
+\cdots+z_k(T_\CC M, \xi_\CC)\varepsilon_2^k. \ee Thus \be
P_1(\tau)
 =2^{4k}\left[z_0(T_\CC M, \xi_\CC)(8\delta_1)^{2k}+z_1(T_\CC M, \xi_\CC)(8\delta_1)^{2k-2}\varepsilon_1
+\cdots+z_k(T_\CC M, \xi_\CC)\varepsilon_1^k \right].\ee

By direct computations, we have \be z_0(T_\CC M, \xi_\CC)
=\left\{\widehat{A}(TM,\nabla^{TM})\right\}^{(8k)}, \ee

\be z_1(T_\CC M,
\xi_\CC)=-\left\{\widehat{A}(TM,\nabla^{TM})\left[48k-\mathrm{ch}\left(B_1(T_\mathbf{C}M,
\xi_\CC)\right)\right] \right\}^{(8k)}.\ee

As we did in (5.49), explicitly expanding $\Theta_1(T_\CC
M,\xi_\CC)$, we get \be \Theta_1(T_\CC M,\xi_\CC)=1+[2(T_\CC
M-(8k)-\xi_\CC+2)-(\xi_\CC\otimes\xi_\CC-2\xi_\CC\wedge
\xi_\CC-\CC^2)]q+\cdots, \ee where the ``$\cdots$" terms are the
terms involving $q^{j\over2}$'s with $j \geq 3.$

By (5.35), (5.55) and (5.62), setting $q=0$ in (5.59), we have
([8, 9]) \be
\left\{\frac{\widehat{L}(TM,\nabla^{TM})}{\cosh^2{\left({c\over2}\right)}}\right\}^{(8k)}
=\sum_{r=0}^{k}2^{6k-6r}z_r(T_\CC M, \xi_\CC).\ee By (5.35),
(5.55) and (5.62), comparing the coefficients of $q$ in (5.59), we
have \be
\begin{split}&\left\{
\frac{\widehat{L}(TM,\nabla^{TM})}{\cosh^2{\left({c\over2}\right)}}\mathrm{ch}\left(2(T_\CC
M, \nabla^{T_\CC
M})-2\cdot8k-2(\xi_\CC-2)-(\xi_\CC\otimes\xi_\CC-2\xi_\CC\wedge
\xi_\CC-\CC^2)\right)
  \right\}^{(8k)}\\
  =&\sum_{r=0}^{k}2^{6k-6r}(48k-64r)z_r(T_\CC M,
\xi_\CC).\\
  \end{split}
\ee

Thus similar to (5.38), we have \be
\begin{split}&\left\{
\frac{\widehat{L}(TM,\nabla^{TM})}{\cosh^2{\left({c\over2}\right)}}\mathrm{ch}\left((T_\CC
M, \nabla^{T_\CC
M})-8k-(\xi_\CC-2)-{1\over2}(\xi_\CC\otimes\xi_\CC-2\xi_\CC\wedge
\xi_\CC-\CC^2)\right)
  \right\}^{(8k)}\\
  =&\ -8k\left(\sum_{r=0}^{k}2^{6k-6r}z_r(T_\CC M, \xi_\CC)\right)+
32\cdot2^6\sum_{r=0}^{k-1}(k-r)2^{6k-6r-6}z_r(T_\CC M, \xi_\CC).
  \end{split}
\ee

Combining (5.63), (5.65) and (5.53), one has \be
\left\{\frac{\widehat{L}(TM, \nabla^{TM})\left[\mathrm{ch}(T_\CC M,
\nabla^{T_\CC
M})-\sinh^2{({c\over2})}\mathrm{ch}\left(2\xi_\CC\oplus\CC^8\right)\right]}{\cosh^2{({c\over2})}}\right\}^{(8k)}\ee
$$=2^{11}[\sum_{r=0}^{k-1}(k-r)2^{6(k-r-1)}z_r(T_\CC M, \xi_\CC)], $$
Which is just (2.32).

\begin{remark} Our main results are obtained by comparing
coefficients of $q$ and $q^2$ in (5.20), (5.31) and coefficients
of $q$ in (5.48), (5.59). It is interesting to examine other
coefficients of higher power of $q$ to get further divisibility
and congruence results. These will be developed elsewhere.

\end{remark}
\section {Acknowledgments} We are grateful to
Professor Weiping Zhang for helpful suggestions and inspiring
discussions with us. Professor M. Atiyah and Pofessor F.
Hirzebruch are deeply appreciated for communications with us. We
also thank Professor Peter Teichner and Professor Nicolai
Reshetikhin for their interests, encouragements and many helpful
discussions.

\bibliographystyle{amsplain}

\begin{thebibliography}{10}


\bibitem {A} L. Alvarez-Gaum\'e and E. Witten, Gravitational
anomalies. {\it Nucl. Phys.} B234 (1983), 269-330.

\bibitem {A} M. F. Atiyah, $K-theory$. Benjamin, New York, 1967.

\bibitem {A} M. F. Atiyah and F. Hirzebruch, Riemann-Roch
theorems for differentiable manifolds. {\it Bull. Amer. Math.
Soc.} 65 (1959), 276-281.

\bibitem {A} M.F. Atiyah and I.M. Singer, The index of elliptic
operators, III, {\it Ann. Math.} 87 (1968), 546-604.

\bibitem {C} K. Chandrasekharan, {\it Elliptic Functions}. Springer-Verlag,
1985.

\bibitem {F} S. M. Finashin, A Pin$^-$-cobordism invariant and a
generalization of Rokhlin signature congruence. {\it Leningrad
Math. J.} 2 (1991), 917-924.

\bibitem {G} Peter B. Gilkey, {\it Invariance Theory, the Heat Equation and the Atiyah-Singer Index Theorem}, Second Edition.
CRC Press, Inc, 1995.

\bibitem {H} F. Han and W. Zhang, Spin$^{c}$-manifold and elliptic genera. {\it
C. R. Acad. Sci. Paris, S$\acute{e}$rie I.} 336 (2003), 1011-1014.
\bibitem {H} F. Han and W. Zhang, Modular invariance, characteristic numbers
and $\eta$ invariants. {\it Journal of Differential Geometry.} 67
(2004), 257-288.

\bibitem {H} F. Hirzebruch, T. Berger and R. Jung, {\it Manifolds and Modular Forms.}
Aspects of Mathematics, vol. E20, Vieweg, Braunschweig, 1992.

\bibitem {H} F. Hirzebruch, {\it Topological Methods in
Algebraic Geometry.} Springer-Verlag, 1966.

\bibitem {H} F. Hirzebruch, Mannigfaltigkeiten
und Modulformen. {\it Jahresberichte der Deutschen Mathematiker
Vereinigung}, Jber. d. Dt. Math.-Verein, 1992, pp. 20-38.

\bibitem {H} Boyuan Hou and Boyu Hou, {\it Differential Geomerty for Physicists}, Second Edition (in Chinese).
Science Press, China, 2004.

\bibitem {L} P. S. Landweber, Elliptic cohomology and modular forms. in {\it
Elliptic Curves and Modular Forms in Algebraic Topology, } p.
55-68. Ed. P. S. Landweber. Lecture Notes in Mathematics Vol.
1326, Springer-Verlag (1988).

\bibitem {L} G. Laures, $K(1)$-local topological modular forms.
{\it Invent Math}. (2004), 371-403.

\bibitem {L} K. Liu, Modular invariance and characteristic
numbers. {\it Commun. Math. Phys}. 174 (1995), 29-42.

\bibitem {L} K. Liu, On Modular Invariance and Rigidity Theorems,
{\it Ph.D Dissertation at Harvard University}. 1993


\bibitem {L} K. Liu and W. Zhang, Elliptic genus and
$\eta$-invariants. {\it Inter. Math. Res. Notices} No. 8 (1994),
319-328.



\bibitem {O} S. Ochanine, Signature modulo 16, invariants de
Kervaire g\'eneralis\'es et nombre caract\'eristiques dans la
$K$-th\'eorie reelle. {\it M\'emoire Soc. Math. France}, Tom. 109
(1987), 1-141.

\bibitem {Z} W. Zhang, Spin$^c$-manifolds and Rokhlin
congruences.
 {\it   C. R. Acad. Sci. Paris,
S\'erie I}, 317 (1993), 689-692.

\bibitem {Z} W. Zhang, Circle bundles, adiabatic limits of
$\eta$ invariants and Rokhlin cogruences. {\it Ann. Inst. Fourier
} 44 (1994), 249-270.

\bibitem {Z} W. Zhang, {\it Lectures on Chern-Weil Theory and
Witten Deformations.} Nankai Tracts in Mathematics Vol. 4, World
Scientific, Singapore, 2001.
$$\  $$

\end{thebibliography}

\end{document}